\documentclass[reqno,12pt,letterpaper]{amsart}
\usepackage{amsmath,amssymb,amstext,amsthm,graphicx,mathrsfs,url,bbm,MnSymbol}
\usepackage[dvipsnames]{xcolor}

\usepackage{tikz-cd}

\usepackage[colorlinks=true,linkcolor=Red,citecolor=Green]{hyperref}

\def\?[#1]{\textbf{[#1]}\marginpar{\Large{\textbf{??}}}}
\newtheorem{thm}{Theorem}
\newtheorem{prop}{Proposition}

\newtheorem{lem}[prop]{Lemma}

\newtheorem{rem}[prop]{Remark}

\numberwithin{equation}{section}
\numberwithin{prop}{section}

\renewcommand{\Re}{\mathop{\rm Re}\nolimits}

\DeclareMathOperator{\WF}{WF}
\DeclareMathOperator{\Res}{{\rm Res}}

\newcommand{\la}{{\lambda}}

\newcommand{\RR}{{\mathbb R}}
\newcommand{\CC}{{\mathbb C}}

\newcommand{\ZZ}{{\mathbb Z}}

\newcommand{\eps}{\varepsilon}

\begin{document}
\title{Twisted Poincar\'e series at zero on surfaces}

\author{Tristan Humbert}
\email{humbertt@imj-prg.fr}
\address{Institut de mathématiques de Paris rive gauche, Sorbonne Université, Campus Pierre et Marie Curie, 4 place Jussieu, 75005, Paris, France.}
\author{Thibault Lefeuvre}
\email{thibault.lefeuvre1@universite-paris-saclay.fr}
\address{Bâtiment 307, rue Michel Magat
Faculté des Sciences d’Orsay, Université Paris-Saclay 
F-91405 Orsay Cedex
}

\author{Sebastián Muñoz-Thon}
\email{sebastian.munoz-thon@universite-paris-saclay.fr}
\address{Bâtiment 307, rue Michel Magat
Faculté des Sciences d’Orsay, Université Paris-Saclay 
F-91405 Orsay Cedex
}
\author{Zhongkai Tao}
\email{ztao@ihes.fr}
\address{Institut des Hautes \'Etudes Scientifiques, 35 route de Ghartres, 91440 Bures-sur-Yvette, France}%

\begin{abstract}
Let $(\Sigma,g)$ be a closed, connected, oriented, Riemannian Anosov surface. We show that there exists an open subset $\mathcal V_g$ of the space of finite-dimensional irreducible
representations of the fundamental group of $\Sigma$, whose complement
has complex codimension at least one, such that the following holds: for all $x,y\in \Sigma$, $\rho \in \mathcal{V}_g$, the meromorphic extension of the twisted Poincaré series $s \mapsto \eta_{g,\rho,x,y}(s)$ vanishes at $s=0$. We also give examples of non-trivial representations whose Poincaré series does not vanish at zero. In addition, we describe the structure of twisted resonant states at zero for these examples and show that they admit non-trivial Jordan blocks.
\end{abstract}

\maketitle

\section{Introduction}

\subsection{Statement of the results}
Let $(\Sigma,g)$ be a closed, connected, oriented, Riemannian surface. Let $\varphi_t \colon S\Sigma\to S\Sigma$ be the geodesic flow on the unit tangent bundle $S\Sigma$ of $\Sigma$. Throughout this article, we will further assume that $(\varphi_t)_{t \in \RR}$ is \emph{Anosov} (e.g., $g$ has negative sectional curvature). We say that $g$ is an \emph{Anosov metric} or that $(\Sigma,g)$ is an Anosov surface in this case.

\subsubsection{Poincaré series} Given $r\in \mathbb{Z}_{\geq 0}$, we introduce $\mathcal{M}_r := \mathrm{Hom}(\pi_1(\Sigma), \mathrm{GL}_r(\CC))$, and let $\mathcal{M}_r^{\mathrm{irr}}$ denote the subset of irreducible representations. For $\rho \in \mathcal{M}_r$, the associated \emph{flat vector bundle} is 
\begin{equation}
    \label{eq:flateq}
     E_\rho:=\widetilde{\Sigma}\times \mathbb C^r /\sim ,\quad (x,\xi)\sim (\gamma(x),\rho(\gamma)\xi), \ \forall \gamma \in \pi_1(\Sigma),
\end{equation}
where $\widetilde{\Sigma}$ is the universal cover of $\Sigma$. We let $\nabla_\rho$ and $d^{\nabla_\rho}$ denote respectively the twisted flat connection and the twisted exterior derivative. We emphasize that we use the convention that the holonomy along $\gamma$ is $\rho(\gamma)^{-1}$.

Let $x,y\in \Sigma$ and $\mathcal P_{x,y}$ be the set of $g$-geodesics starting at $x$ and ending at $y$. When $x=y$, the set $\mathcal{P}_{x,x}$ also contains the trivial loop reduced to the singleton $\{x\}$. For $\gamma \in \mathcal P_{x,y}$, the \emph{parallel transport} along $\gamma$ with respect to $\nabla_\rho$ is a linear map
\[
\tau_\rho(\gamma)\in \mathrm{Hom}(E_\rho(x), E_\rho(y)). 
\]
We define the \emph{twisted Poincaré series} for $s\in \CC$, $\mathrm{Re}(s)\gg 1$, by
\begin{equation}
\label{eq:Poinc}
\eta_{g,\rho,x,y}(s):=\sum_{\gamma\in \mathcal P_{x,y}}e^{-s\ell_g(\gamma)}\tau_\rho(\gamma)\in \mathrm{Hom}(E_\rho(x), E_\rho(y)),
\end{equation}
where $\ell_g(\gamma)$ is the length of $\gamma$. The twisted Poincaré series \eqref{eq:Poinc} converges and is holomorphic for $\Re(s) \gg 1$. It admits a meromorphic extension to $\CC$ (Proposition \ref{prop:mero}). 

\subsubsection{Generic value at $0$.}

In this article, we study the value of the meromorphic extension of $\eta_{g,\rho,x,y}(s)$ at $s=0$. Our first main result is the following:
\begin{thm}[Generic case]
\label{theo:1}
Let $(\Sigma,g)$ be a closed, connected, oriented, Anosov surface. There exists an open subset $\mathcal V_g \subset \mathcal{M}_r^{\mathrm{irr}}$, invariant by conjugation, containing all non-trivial unitary representations, whose complement is of complex dimension $\geq 1$, such that the following holds: for all $\rho \in \mathcal V_g$, $x, y \in \Sigma$,
    \[
    \eta_{g,\rho,x,y}(0)=0.
    \]
\end{thm} 

For $r=1$ and $\rho=\rho_{\mathrm{triv}}$, the trivial representation, Dang and Rivière showed \cite[Theorem 1.2]{DR24} that $\eta_{g,\rho_{\mathrm{triv}},x,y}(0)=1/\chi(\Sigma)$, where $\chi(\Sigma)$ is the Euler characteristic of $\Sigma$\footnote{When $x= y$, they obtain $\eta_{g,\rho_{\mathrm{triv}},x,x}(0)=\chi(\Sigma)^{-1}-1$. The discrepancy comes from the fact that they do \emph{not} include the trivial loop in their Poincaré series when the two points are equal.}. Theorem \ref{theo:1} shows that $\rho \mapsto \eta_{g,\rho,x,y}(0)$ is \emph{not} continuous with respect to the representation $\rho\in \mathcal{M}_r^{\mathrm{irr}}$. The proof uses the work of the first and last-named authors \cite{HT26} which describes the \emph{twisted Pollicott--Ruelle resonant spaces at zero} for a generic representation.

\subsubsection{Holomorphic rank one connections} Let $K \to \Sigma$ be the canonical line bundle of $\Sigma$, seen as complex curve, where the complex structure is the one induced by the conformal class of the metric $g$, see \cite[Chapter 15]{Lef} for details. Recall that a holomorphic $1$-form $\omega \in H^0(\Sigma; K) $ defines a rank one representation $\rho_\omega \in \mathcal{M}_1$ by
\begin{equation}
    \label{eq:rhow}
    \rho_\omega(\gamma):= \exp\left(-\int_\gamma \omega\right), \qquad \forall \gamma \in \pi_1(\Sigma).
\end{equation}
The corresponding flat complex line bundle $E_{\rho_\omega} \to \Sigma$ can be identified with the trivial line bundle $\CC \times \Sigma$ equipped with the connection $d+\omega \wedge$.

For any $x\in \Sigma$, let
\begin{equation}
    \label{eq:ev_x}
\mathrm{ev}_x^\omega:H^0(\Sigma;E_{\rho_\omega})\to E_{\rho_\omega}(x),\quad s\mapsto s(x),
\end{equation}
be the evaluation map, where $H^0(\Sigma;E_{\rho_\omega})$ denotes the holomorphic sections of the line bundle $E_{\rho_\omega}$. Note that $H^0(\Sigma;E_{\rho_\omega}) \simeq \CC$ and \eqref{eq:ev_x} is an isomorphism (Lemma \ref{lemma:evaluation}). In particular, this defines a canonical section
\[
\sigma_\omega(x,y) := \mathrm{ev}_y^\omega \circ (\mathrm{ev}_x^\omega)^{-1} \in \mathrm{Hom}(E_{\rho_\omega}(x), E_{\rho_\omega}(y)).
\]
The following holds:

\begin{thm}[Holomorphic 1-forms]
    \label{theo:2}
    Let $(\Sigma,g)$ be a hyperbolic surface. There exists a neighborhood $U$ of $0$ in $H^0(\Sigma; K)$ such that for any $x,y\in \Sigma$ and $\omega \in U$:
    \begin{equation}
\label{eq:counterexintro}
\eta_{g,\rho_{\omega},x,y}(0)= \sigma_\omega(x,y)/\chi(\Sigma).
\end{equation}
\end{thm}

The equality \eqref{eq:counterexintro} provides a nonzero value at $s=0$ and contrasts with Theorem \ref{theo:1}. We also note that the space of holomorphic $1$-forms on $\Sigma$ has complex dimension $\textsl{g}$, where $\textsl{g}\geq 2$ is the genus of $\Sigma$, while the space of rank-one representations has complex dimension $2\textsl{g}$. Thus, our set of examples is contained in a subset of complex codimension $\textsl{g}\geq 2$, which is consistent with the results of Theorem \ref{theo:1}. Similarly to Theorem \ref{theo:2}, there should also exist higher rank (non-trivial) representations ($r \geq 2$) for which $\eta_{g,\rho,x,y}(0) \neq 0$.

\subsubsection{Jordan blocks} The overall strategy of the proof builds on \cite{DR24} and is based on a description of the \emph{twisted Pollicott--Ruelle resonances} at zero, see \S \ref{sec:Strategy} for details. Our last result describes the resonant spaces for the examples constructed in Theorem \ref{theo:2} and shows that it differs from the generic picture obtained in \cite[Theorem 2]{HT26}. In the next theorem, we let $\Res_0^{k,\ell}(\rho,0) $ denote the space of generalized resonant states at $0$ in the kernel of $\mathcal L_{X^{\rho}}^\ell$, for the action of $\mathcal L_{X^{\rho}}$ on $k$-forms, see \S \ref{sec:Ruelle} for more details.
\begin{thm}
    \label{theo:3}
    Let $(\Sigma,g)$ be a hyperbolic surface of genus $\textsl{g}\geq 2$. There exists a neighborhood $U$ of $0$ in $H^0(\Sigma; K)$, such that for any $\omega \in U\setminus \{0\}$:
 \begin{align*}&\mathrm{dim}(\Res_0^{k,\infty}(\rho_{\omega},0))=1, \quad k=0,2, 
\\ &\mathrm{dim}(\Res_0^{1,\infty}(\rho_{\omega},0))=2 \textsl{g},\quad \mathrm{dim}(\Res_0^{1,1}(\rho_{\omega},0))=2 \textsl{g}-1. 
\end{align*}
\end{thm}
To the best of our knowledge, Theorem \ref{theo:3} provides the first examples of representations of $\pi_1(\Sigma)$ that admit a Jordan block at zero. Let us note that examples of representations of $\pi_1(S\Sigma)$ which do \emph{not} factor through $\pi_1(\Sigma)$ admitting a Jordan block at zero were constructed in \cite[Theorem 4]{HT26}.

\subsection{Related works}
The value at $s=0$ of the Poincaré series for Anosov surfaces was computed by Dang and Rivière in \cite[Theorem 1.2]{DR24}. For hyperbolic surfaces, an integrated formula over all points $x,y\in \Sigma$ was obtained before by Paternain \cite{Pat00}. The result of \cite{DR24} was extended by Chaubet \cite[Theorem 2]{Cha22} to the setting of negatively curved surfaces with totally geodesic boundary. In the same paper, Chaubet studies the Poincaré series counting orthogeodesics joining boundary points and shows that it vanishes at $s=0,$ see \cite[Theorem 1]{Cha22}. See also \cite{CGB24} for a more general statement relating values of some dynamical
series at $s=0$ to some topological invariants and \cite{DLR23,DGBLR25} for results on Poincaré series for convex bodies. Finally, the Poincaré series associated to Anosov representations was shown to have a meromorphic extension to $\mathfrak{a}^*_{\CC}$ by Guedes Bonthonneau, Weich and the second-named author, see \cite[Theorem 1.10]{GBLW26}.

\subsection{Strategy of the proof}
\label{sec:Strategy}
We follow the general strategy of \cite{DR24} and express the Poincaré series as a distributional pairing involving the resolvent of the geodesic vector field, see Proposition \ref{prop:DR}. We then deduce that the Poincaré series admits a meromorphic extension to $\mathbb C$ from the theory of (twisted) \emph{Pollicott-Ruelle resonances}--this is the content of Proposition \ref{prop:mero}. The behavior of the meromorphic extension at $s=0$ is determined by the \emph{resonant spaces at zero for the action on forms of degree $k=0,1,2$ in the kernel of the contraction operator $\iota_X$}.

In the non-twisted case, Dang and Rivière used the work of Dyatlov and Zworski \cite{zazi}, who showed that there are no Jordan blocks at $0$ and that resonant $1$-forms are closed. Chaubet used in \cite{Cha22} the results of Hadfield \cite{Had18}, which extends the work of Dyatlov and Zworski to the setting of surfaces with boundary. The description of the resonant states at zero for a generic representation of $\pi_1(\Sigma)$ was recently obtained by the first and last named authors in \cite{HT26} and we briefly recall in \S \ref{sec:Ruelle} these results. The last step consists in computing the value at $s=0$ by constructing explicit primitives of the distributional currents that appear in the formula. Using that $\rho$ is nontrivial and belongs to the open set described in Proposition \ref{prop:HT}, we obtain Theorem \ref{theo:1}. In Proposition \ref{prop:gen-res-states}, we give an explicit description of the Jordan block structure at zero using an idea coming from \emph{quantum-classical correspondence} for hyperbolic surfaces \cite{GHW}.

\subsection*{Acknowledgement}
We thank Gabriel Rivière for discussions related to this project. T.H., T.L. and S. M.-T. were supported by the European Research Council (ERC) under the European Union’s Horizon 2020 research and innovation programme (Grant agreement no. 101162990 — ADG). Z.T.~was partially supported by a grant from the National Research Agency (ANR) under the France 2030 program, reference ANR-25-CFFS-0012.

\subsection*{AI disclosure}
Theorem \ref{theo:1} was proved without generative AI. We then asked an LLM to search for examples of representations with non-vanishing Poincaré series at zero. More precisely, we asked the LLM to focus on 1-dimensional representations, and obtain an example for hyperbolic metrics using a twisted Selberg trace formula. The idea to consider holomorphic 1-forms in the statement of Theorem~\ref{theo:2} is due to the LLM: it produced examples using explicit computations involving the twisted Selberg trace formula. We then found a shorter proof of Theorem~\ref{theo:2}, which is the one presented in the article, using Pollicott--Ruelle resonances. This led us to prove Theorem \ref{theo:3} on the precise Jordan block structure of the examples. We emphasize that this last result does \emph{not} follow from the original argument of the LLM. All proofs were written by the authors. Generative AI was not used in the writing, nor in the proofreading of this manuscript.

\section{Preliminaries}

\subsection{Unit tangent bundle} \label{subsection:unit-tangent-bundle}
Let $(\Sigma,g)$ be a closed oriented Riemannian Anosov surface.

\subsubsection{Dynamical splittings} Let $M := S\Sigma$ be its unit tangent bundle, and $\pi \colon M \to \Sigma$ be the footpoint projection. Let $X$ be the geodesic vector field. The Anosov property means that there exist a continuous flow-invariant splitting $TM=E_u\oplus \mathbb R X\oplus E_s$ of the tangent bundle and constants $C,\vartheta>0$ such that: 
\[
\|d\varphi_t|_{E_s}\|\leq Ce^{-\vartheta t},\, t\geq 0, \qquad \|d\varphi_t|_{E_u}\|\leq Ce^{\vartheta t},\, t\leq 0.
\]
The bundle $E_s$ (resp. $E_u$) is the \emph{stable} (resp. \emph{unstable}) bundle of $(\varphi_t)_{t\in \mathbb R}.$ The Anosov splitting induces a dual splitting on the cotangent bundle of $M$:
\begin{equation}
\label{eq:dualdecomp}
T^*M=E_u^*\oplus \mathbb R \alpha \oplus E_s^*,\quad E_u^*(E_u\oplus \mathbb RX)=0, \ E_s^*(E_s\oplus \mathbb RX)=0,
\end{equation}
where $\alpha\in C^{\infty}(M;T^*M)$ is the contact $1$-form on $M$, which satisfies $\iota_X\alpha  = 1$, $\iota_Xd\alpha=0$.

\subsubsection{Vertical and horizontal bundles} An oriented Riemannian surface admits an integrable almost-complex structure, i.e., there is a parallel section $J\in C^\infty(M; \mathrm{End}(TM))$ such that $J^2=-\mathrm{Id}$, and the area form associated to $g$ is given by $dA_g=g(J\bullet,\bullet)$. 
One defines the fiberwise rotation of angle $\theta$ by
$$\rho_{\theta}\colon SM\to SM,\quad \rho_\theta(x,v)=(x,e^{\theta J}v), $$
where $e^{\theta J}v$ is the unit vector obtained by rotating $v$ by an angle $\theta$ in the positive direction (with respect to the orientation of $M$).
The \emph{vertical vector field} $V$ is defined as the generator of $\rho_\theta$. The \emph{horizontal vector field} is then defined by $H:=[V,X]$. We let $\mathbb V:=\mathbb R V$ (resp. $\mathbb H:=\mathbb R H$) be the vertical bundle (resp. horizontal bundle). Observe that $TM=\mathbb R X\oplus \mathbb V\oplus \mathbb H$. We will use the corresponding decomposition of the cotangent bundle:
\begin{equation}
\label{eq:VHsplit}T^*M=\mathbb H^*\oplus \mathbb R \alpha \oplus \mathbb V^*,\quad \mathbb V^*(\mathbb H\oplus \mathbb RX)=0, \quad \mathbb H^*(\mathbb V\oplus \mathbb RX)=0.
\end{equation}

\subsubsection{Fourier decomposition} \label{sssection:fourier}

Let 
$$\Omega_k:=\{u\in C^{\infty}(S\Sigma): Vu=iku\},\quad k\in \mathbb Z. $$
Each smooth function $f\in C^{\infty}(S\Sigma)$ can be decomposed uniquely into a sum of Fourier modes:
\[
f = \sum_{k \in \ZZ} f_k, \qquad f_k \in \Omega_k.
\]
See \cite[Equation (15.2.5)]{Lef}. The \emph{Guillemin--Kazhdan raising and lowering operators} are defined by $\eta_\pm := (X\mp iH)/2$ and map
\[
\eta_\pm \colon \Omega_k \to \Omega_{k \pm 1}, \qquad k \in \ZZ.
\]

\subsubsection{Contact distribution} The contact distribution is the flow-invariant subbundle
\begin{equation}
\label{eq:contact}\mathbb H\oplus \mathbb V=\mathrm{ker}(\alpha)=E_u\oplus E_s.
\end{equation}
Since the flow is Anosov, the vertical direction is transverse to both the stable and unstable directions everywhere (see \cite[Theorem 13.3.6]{Lef}). As before, we will use the dual statement:
\begin{equation}
    \label{eq:Vtransverse}
    (\mathbb H^*\oplus \mathbb R\alpha) \cap E_u^*=(\mathbb H^*\oplus \mathbb R\alpha) \cap E_s^*=\{0\}.
\end{equation}
Next, we use the fact that Anosov flows have no conjugate points (\cite[Theorem 13.3.10]{Lef}). In the cotangent bundle, this statement reads:
\begin{equation}
    \label{eq:noconj}
    \forall t\neq 0,\quad d\varphi_t^{\top}(\mathbb H^*)\cap \mathbb H^*=\{0\}.
\end{equation}
From now on, we orient $M$ using the contact volume form $\alpha \wedge d\alpha$. With this choice of orientation,
\begin{equation}
    \label{eq:orientation}
    (X,V,H) \text{ is a positively orthonormal basis of } TM,
\end{equation}
with respect to the Sasaki metric. The following fact will be used in the proof of Theorem~\ref{theo:1}:

\begin{lem}
There exists $a_t,b_t\in C^\infty(M)$ such that $a_t > 0$ $($resp. $a_t < 0)$ for $t > 0$ $($resp. $t < 0)$ and 
\begin{equation}
    \label{eq:positive}
     d\varphi_t(V)=a_tH+b_tV.
\end{equation} 
\end{lem}

\begin{proof}We first use the flow-invariance of the contact distribution \eqref{eq:contact} to deduce that there exists $a_t,b_t\in C^\infty(M)$ such that \eqref{eq:positive} holds. In fact, the coefficients are given by \emph{Jacobi fields} and the absence of conjugate points recalled above implies that $a_t$ cannot vanish for $t \neq 0$, see for instance \cite[Lemma 13.2.2]{Lef}. We only need to determine the sign of $a_t$. For this, we Taylor expand near $t=0$: $d\varphi_t(V)=V+t[V,X]+O(t^2)=V+tH+O(t^2)$. This shows that $a_t\sim t$ and completes the proof.
\end{proof}

\subsubsection{Twisted geometry on hyperbolic surfaces}

\label{sssection:twisted-geometry}

We now record useful identities which will be used in the proofs of Theorem \ref{theo:2} and Theorem \ref{theo:3}. In this paragraph, we further assume that $g$ is hyperbolic. The unstable horocyclic vector field is defined by
\begin{equation}
    \label{equation:horocyclic}
    U_- := H + V,
\end{equation}
and spans $E_u$ on $S\Sigma$. It satisfies $[X,U_-]=-U_-$.

Let $\omega$ be a holomorphic $1$-form on $\Sigma$. We will identify $\omega$ with the function on $M=S\Sigma$ (still denoted by $\omega$) given by $(x,v) \mapsto \omega_x(v)$. Define
\[
X^\omega := X + \omega, \qquad H^\omega := H + i \omega, \qquad U_-^\omega := U_- + i\omega, \qquad \eta_+^\omega := \eta_+ + \omega
\]
Finally, the \emph{twisted Casimir operator} is
\begin{equation}
     \label{eq:Casimir}
         C^\omega:= (X^\omega)^2 + (H^\omega)^2-V^2 = (X^\omega)^2-X^\omega+(U_-^\omega)^2-2VU_-^\omega.
     \end{equation}
The following holds:

\begin{lem}
\label{lemma:casimir} We have
$[X^\omega,U^{\omega}_-] = - U^\omega_-$ and $[C^\omega,V]=0$.
\end{lem}

\begin{proof}
    Immediate using the commutation relations $[X,H]=-V$, $[V,X]=H$, $[H,V]=X$ (see \cite[Lemma 15.2.1]{Lef}).
\end{proof}

\subsection{Twisted Pollicott--Ruelle resonances}
\label{sec:Ruelle}

We refer to \cite[\S 2.1]{HT26} for a more detailed exposition on twisted Pollicott--Ruelle resonances. Let $\rho \in \mathcal{M}_r$ and let $\mathcal E_{\rho} := \pi^*E_\rho$ denote the pullback of the associated flat vector bundle over $M=S\Sigma$ with differential $d^{\nabla_\rho}$, see \eqref{eq:flateq}. Define the graded vector bundle
\begin{equation}
\label{eq:forms}\mathscr E_\rho:=\bigoplus_{k=0}^{3} \mathscr E^k_\rho,\quad \mathscr E^k_\rho=\Lambda^kT^*M\otimes \mathcal E_\rho, \quad \mathscr E^k_{\rho,0}:=\mathscr E^k_\rho\cap \mathrm{ker}(\iota_X). 
\end{equation}

\subsubsection{Resolvent} For $k=0,1,2$, the generator $X$ acts on sections of the smooth vector bundle $C^{\infty}(M; \mathscr E^k_{\rho,0})$ by
\[
\mathcal L_{X^\rho}:=\iota_X d^{\nabla_\rho}+d^{\nabla_\rho} \iota_X.
\]
It follows from \cite{Fau10} that the resolvent defined for $\Re(s) \gg 1$ by
\begin{equation} \label{eq:resolvent-integral}
    R_{k}^{\rho}(\lambda):=(\mathcal L_{X^\rho}+\lambda)^{-1}=\int_{0}^{\infty}e^{-t(\mathcal{L}_{X^{\rho}}+\la)}dt,
\end{equation}
admits a meromorphic extension to $\mathbb C$ as a map:
\begin{equation}
\label{eq:Rk}
    R_k^\rho(\lambda) \colon C^{\infty}(M; \mathscr E^k_{\rho,0})\to \mathcal D'(M; \mathscr E^k_{\rho,0}). 
    \end{equation}
The poles are intrinsic and called \emph{twisted Pollicott--Ruelle resonances} (on $\mathscr E^k_{\rho,0}$). The set of resonances is denoted by $\mathrm{Res}^k(X,\rho)$ or $\mathrm{Res}^k(\rho)$ when there is no confusion on the metric.

\subsubsection{Resonant states} A complex number $\lambda_0\in \mathbb C$ is in $ \mathrm{Res}^k(\rho)$, if and only if
\begin{equation}
\label{eq:res}
    \Res_0^{k,1}(\rho,\lambda_0):=\{u\in \mathcal{D}'(M;\mathscr E^k_{\rho,0}): \mathcal{L}_{X^\rho} u = -\lambda_0 u , \, \WF(u)\subset E_u^* \}\neq \{0\},
\end{equation}
where $\mathrm{WF}(u)$ denotes the wavefront set (see \cite[Chapter 3]{Hor}). The non-zero elements of $\Res_0^{k,1}(\rho,\lambda_0)$ are called \emph{resonant states} at $\lambda_0$. The spectral projector at $\lambda_0$ is given by
\begin{equation}
\label{eq:Pi_k}\Pi^\rho_k(\lambda_0)=\frac{1}{2i\pi}\int_{\gamma}R_k^\rho(z)dz,
\end{equation}
where $\gamma$ is a small loop around $\lambda_0$. The (algebraic) multiplicity of $\lambda_0$ is given by the rank of the spectral projector.  The \emph{generalized resonant states} are the elements in the range of $\Pi^k_\rho(\lambda_0)$:
\begin{equation}
\label{eq:genres}
\Res_0^{k,\infty}(\rho,\lambda_0)=\{u\in \mathcal{D}'(M;\mathscr E^k_{\rho,0}): \exists \ell\in \mathbb{Z}_{\geq 0},\ (\mathcal{L}_{X^\rho}+\lambda_0)^{\ell} u = 0 , \, \WF(u)\subset E_u^* \}.
\end{equation}
We will say that the resonance $\lambda_0$ has \emph{no Jordan block} if $\Res_0^{k,\infty}(\rho,\lambda_0)=\Res_0^{k,1}(\rho,\lambda_0)$.

The resolvent $R_k^{\rho}(\lambda)$ and the spectral projection $\Pi_k^{\rho}$ are also defined on $C^{\infty}(M;\mathscr E_{\rho}^k)\to \mathcal D'(M;\mathscr E_{\rho}^k)$, and the space of corresponding resonant states is
\begin{equation}\label{eq:res-1}
    \Res^{k,1}(\rho,\lambda_0) = \{u\in \mathcal{D}'(M;\mathscr E^k_{\rho}): \mathcal{L}_{X^\rho} u =- \lambda_0 u , \, \WF(u)\subset E_u^* \}.
\end{equation}
Note that we have the following relation between ranges of the two spectral projectors:
\begin{equation}
    \label{eq:relalpha}
    \forall k=0,1,2,3\quad \Res^{k,1}(\rho,\lambda_0)=\Res^{k,1}(\rho,\lambda_0)\wedge \alpha\oplus \Res^{k,1}_0(\rho,\lambda_0),
\end{equation}
where $\Res^{-1,1}_0(\rho,\lambda_0)=0$ by convention. 
We will need a result due to the first and last named authors \cite[Theorem 2]{HT26}.

\begin{prop}
\label{prop:HT}
Let $(\Sigma,g)$ be a closed, connected, oriented Anosov surface. Then, there exists an open subset $\mathcal U_g \subset \mathcal{M}_r^{\mathrm{irr}}$, invariant by conjugation, whose complement in $\mathcal{M}_r^{\mathrm{irr}}$ has complex codimension $\geq 1$, such that the following holds: for any $\rho \in \mathcal{U}_g$:
    \begin{equation}
    \label{eq:dim}
    \begin{split}
        &\mathrm{dim}(\Res_0^{k,\infty}(\rho,0))=0, \quad k=0,2, 
        \\ &\mathrm{dim}(\Res_0^{1,\infty}(\rho,0))=\mathrm{dim}(\Res_0^{1,1}(\rho,0))=-\mathrm{dim}(\rho)\chi(\Sigma). 
        \end{split}
          \end{equation}
      Moreover, the resonant states in $\Res_0^{1,1}(\rho,0)$ are $d^{\nabla_\rho}$-closed.
\end{prop}

That the resonant states are $d^{\nabla_\rho}$-closed can be deduced from the proof of \cite[Lemma 4.5]{HT26}. Finally, we define
\begin{equation}
    \label{def:Vg}
\mathcal{V}_g \subset \mathcal{M}_r
    \end{equation}
to be the set of all irreducible representations such that \eqref{eq:dim} hold, and the resonant states in $\Res_0^{1,1}(\rho)$ are closed. By the proof of \cite[Lemma 4.5]{HT26}, we know $\mathcal{V}_g$ is open. By Proposition~\ref{prop:HT} and  \cite[Proposition 4.6]{HT26}, $\mathcal{U}_g\subset \mathcal{V}_g$, and $\mathcal{V}_g$ contains all nontrivial irreducible unitary representations.

\section{Expressing the Poincaré series in terms of the resolvent}
\label{section:poincare-pollicott-ruelle}

We follow the general strategy of \cite{DR24} and write the twisted Poincaré series in terms of the resolvent \eqref{eq:resolvent-integral}. This will show that the Poincaré series extends meromorphically to $\mathbb C$ and express the value at $s=0$ using the description of the twisted Pollicott--Ruelle resonances at zero obtained in \cite{HT26}.

\subsection{Expressing the twisted Poincaré series} Let $x,y\in \Sigma$,  $\xi \in E_\rho(x)$, $\mu\in E_\rho^*(y)$. We extend these vectors to smooth sections $\xi_\star \in C^\infty(\Sigma; E_\rho)$ (resp. $\mu_\star \in C^\infty(\Sigma; E_\rho)$) such that $\xi_\star(x)=\xi$ (resp. $\mu_\star(y)=\mu$). In the next statement, pairings on the left-hand side are understood in the sense of pairings between vectors and covectors in finite-dimensional vector spaces, while pairings on the right hand-side are in the sense of currents (see \cite[Chapter 10]{Lef}).

\begin{prop}
    \label{prop:DR}
    Let $(\Sigma,g)$ be a closed, connected, oriented, Anosov surface, and $\rho \in \mathcal{M}_r$. Let $x,y\in \Sigma$ and let $\xi \in E_\rho(x)$ and $\mu\in E_\rho^*(y)$. Then, for any $\varepsilon>0$ small enough, for any $s\in \mathbb C$ with $\mathrm{Re}(s)\gg 1$:
    \begin{enumerate}
    \item If $x\neq y$,
    \begin{equation}
\label{eq:pairing}
\langle \eta_{g,\rho,x,y}(s)\xi,\mu\rangle =e^{-s\varepsilon}\langle R_1^\rho(s)\big(\iota_X[S_x\Sigma]\otimes \pi^*\xi_\star\big),e^{\varepsilon \mathcal L_{X^{\rho^*}}}([S_y\Sigma]\otimes \pi^*\mu_\star) \rangle;
\end{equation}
\item If $x=y$, then 
\begin{equation}
\label{eq:pairing2}
\langle \eta_{g,\rho,x,x}(s)\xi,\mu\rangle=e^{-s\varepsilon}\langle R_1^\rho(s)\big(\iota_X[S_x\Sigma]\otimes \pi^*\xi_\star\big),e^{\varepsilon \mathcal L_{X^{\rho^*}}}([S_x\Sigma]\otimes \pi^*\mu_\star) \rangle+\langle \xi,\mu\rangle.
\end{equation}
\end{enumerate}
\end{prop}

\begin{rem}
\label{rem:epsilon}
The number $\varepsilon>0$ is chosen such that $\varphi_{\varepsilon} (S_x\Sigma)\cap S_y\Sigma =\emptyset$. This guarantees that the distributional pairings in the proof are well-defined by wavefront set conditions, see for instance \cite[Chapter 3]{Hor}. In particular, if $ x\neq y$, one can take $\varepsilon=0$ in the previous proposition. The different formulas for $x = y$ and $x \neq y$ are due to the fact that we include the trivial loop in the definition of $\eta_{g,\rho,x,x}(s)$ when $x=y$ $($cf. \cite[Proposition 3.8]{DR24}$)$.
\end{rem}

\begin{proof}
For $s\in \CC$ such that $\mathrm{Re}(s)\gg 1$, the Poincaré series is equal to
\begin{equation}
\label{eq:etauseful}
    \langle \eta_{g,\rho,x,y}(s)\xi,\mu\rangle = \sum_{\gamma\in \mathcal P_{x,y}}e^{-s\ell_g(\gamma)}\langle \tau_\rho(\gamma)\xi,\mu\rangle.\end{equation}
Let $\varepsilon>0$ be small enough, and define $c_{x,\varepsilon}:=\pi(\varphi_\varepsilon(S_x\Sigma))\subset \Sigma$, the geodesic sphere of radius $\varepsilon>0$ centered at $x$. Let $\tau_\rho(\eps) \colon E_\rho(x) \to C^\infty(c_{x,\varepsilon},E_\rho)$ be the map defined by parallel transport of a section $\xi \in E_\rho(x)$ to $c_{x,\varepsilon}$ along radial geodesics. We can define the Poincaré series of $y$ and $c_{x,\varepsilon}$ by:
$$\eta_{g,\rho,c_{x,\varepsilon},y}(s):=\sum_{\gamma \in \mathcal P_{c_{x,\varepsilon},y}}\tau_{\rho}(\gamma)\tau_\rho(\eps)e^{-s\ell_g(\gamma)}, $$
where $\mathcal P_{c_{x,\varepsilon},y}$ denotes the set of \emph{orthogeodesics} joining $c_{x,\varepsilon}$ to $y$, that is, the set of geodesics that are orthogonal to $c_{x,\varepsilon}$, join a point on $c_{x,\varepsilon}$ to $y$, and have outward pointing initial velocity. When $x=y$, we define $\mathcal P_{c_{x,\varepsilon},x}$ to be the set of orthogeodesics of length strictly greater than $\varepsilon$ to ensure that the set of length is discrete (geodesic rays of length $\varepsilon$ joining $x$ to $c_{x,\varepsilon}$ are not included in this set). 

Any non-trivial $\gamma \in \mathcal{P}_{x,y} $ intersects $c_{x,\varepsilon}$ orthogonally by the Gauss lemma. That is, if $\varepsilon>0$ is chosen small enough, $\gamma$ determines a unique orthogeodesic $\gamma_\eps$ from $c_{x,\varepsilon}$ to $y$ whose length is equal to $\ell_g(\gamma)-\varepsilon$. This yields the following formulas for $\varepsilon > 0$ small enough and $s \in \CC$ such that $\Re(s) \gg 1$:

\begin{enumerate}
    \item If $x\neq y$, then
\begin{equation}
    \label{eq:magic0}
    \eta_{g,\rho,x,y}(s)=e^{-s\varepsilon}\eta_{g,\rho,c_{x,\varepsilon},y}(s);
\end{equation}
   \item If $x= y$, then
\begin{equation}
    \label{eq:magic2}
    \eta_{g,\rho,x,x}(s)=e^{-s\varepsilon}\eta_{g,\rho,c_{x,\varepsilon},x}(s)+\mathrm{Id}_{E_\rho(x)},
\end{equation}
\end{enumerate}

By \cite[Lemma 3.6]{DR24}, for any $\eps > 0$ small enough, and $T\geq 0$, if $v \in \varphi_{\varepsilon}(S_x\Sigma) \cap \varphi_{-T}(S_y\Sigma)$, then
\[
T_vM = T_v(\varphi_{\varepsilon}(S_x\Sigma)) \oplus T_v(\varphi_{-T}(S_y\Sigma)) \oplus \RR X(v).
\]
In fact, if $x\neq y$, this holds for $\varepsilon=0$. In particular, $\varphi_{\varepsilon}(S_x\Sigma)$ and $S_y\Sigma$ satisfy the conditions of \cite[Lemma 2.1]{DR24}. Note that $c_{x,\varepsilon}$ is homotopically trivial so it bounds a simply connected open set containing $x$ (namely, the geodesic ball of radius $\varepsilon$ centered at $x$).
Adapting the proof of \cite[Lemma 2.1]{DR24} and using that the propagator $e^{t\mathcal L_{X^\rho}}$ acts on sections of the flat bundle by parallel transport (with respect to $\nabla_\rho$) along the geodesic flowlines, we find that the following equality holds in $\mathcal D'(\mathbb R_+)$, 
\begin{equation}
    \label{eq:DRlemma21}
    \begin{split}
    f_{\xi,\mu,\varepsilon}(t)&:=\sum_{\gamma_\varepsilon \in \mathcal P_{c_{x,\varepsilon},y}}\langle \tau_\rho (\gamma_\varepsilon)\tau_\rho(\eps)  \xi,\mu\rangle\delta(t-\ell_g(\gamma_\varepsilon))
    \\&=\langle e^{-t\mathcal L_{X^\rho}}e^{-\varepsilon\mathcal L_{X^\rho}}(\iota_X[S_x\Sigma]\otimes \pi^*\xi_\star), [S_y\Sigma]\otimes \pi^* \mu_\star)\rangle.
    \end{split}
\end{equation}
By \eqref{eq:orientation} and \eqref{eq:positive}, the term $\epsilon_{\tau}(x)=-1$ in \cite[Equation (7)]{DR24} cancels with $(-1)^{\dim N_1}=-1$.

We pair the above equality with $t\mapsto e^{-st}$ and find\footnote{The pairing is well-defined for $s$ large enough although $t \mapsto e^{-st}$ is not compactly supported.}
\begin{equation}
\label{eq:long}
\begin{split}
    \langle \eta_{g,\rho,c_{x,\varepsilon},y}(s) \xi,\mu\rangle & =\langle  f_{\xi,\mu,\varepsilon}(t),e^{-st}\rangle_{\mathcal{D}'(\mathbb R_+)\times C^\infty(\mathbb R_+)}
    \\& \overset{\eqref{eq:DRlemma21}}{=}\int_{0}^{+\infty}\langle e^{-t\mathcal L_{X^\rho}}e^{-\varepsilon\mathcal L_{X^\rho}}(\iota_X[S_x\Sigma]\otimes \pi^*\xi_\star), [S_y\Sigma]\otimes \pi^*\mu_\star)\rangle e^{-st}dt
    \\&=\left\langle \int_{0}^{+\infty} e^{-t(s+\mathcal L_{X^\rho})} e^{-\varepsilon\mathcal L_{X^\rho}}(\iota_X[S_x\Sigma]\otimes \pi^*\xi_\star) dt, [S_y\Sigma]\otimes \pi^*\mu_\star\right\rangle
    \\& \overset{\eqref{eq:resolvent-integral}}{=}\langle R_1^\rho(s)e^{-\varepsilon \mathcal L_{X^\rho}}\big(\iota_X[S_x\Sigma]\otimes \pi^*\xi_\star\big),[S_y\Sigma]\otimes \pi^*\mu_\star\rangle \\
    & =\langle R_1^\rho(s)\big(\iota_X[S_x\Sigma]\otimes \pi^* \xi_\star\big),e^{\varepsilon \mathcal L_{X^{\rho^*}}}([S_y\Sigma]\otimes \pi^*\mu_\star)\rangle
 \end{split}   
\end{equation}
where we used in the first line the definition of $\eta_{g,\rho,c_{x,\varepsilon},y}(s)$, we exchanged the absolutely converging integrals (for $\Re(s) \gg 1$) in the third line, we used the integral representation of the resolvent when $\mathrm{Re}(s)\gg 1$ in the fourth line, and $[R_1^\rho(s),e^{-\varepsilon \mathcal L_{X^\rho}}]=0$ in the last line. The claimed result then follows from \eqref{eq:magic0} and \eqref{eq:magic2}.
\end{proof}

\begin{prop}
    \label{prop:mero}
    Let $(\Sigma,g)$ be a closed, connected, oriented, Anosov surface. Let $\rho \in \mathcal{M}_r$. Then, the Poincaré series $s \mapsto \eta_{g,\rho,x,y}(s)$ admits a meromorphic extension to $\mathbb C$ with poles contained in $\mathrm{Res}^1(\rho).$
\end{prop}

\begin{proof}
It follows from \cite[Proposition 3.3]{zeta}
that $R_{k}^{\rho}(s) \colon C^{\infty}(M;\mathscr E_{\rho}^k)\to \mathcal{D}' (M;\mathscr E_{\rho}^k)$ has a meromorphic continuation to $\mathbb{C}$, and near every pole $s_0$, 
\begin{equation}
    R_k^{\rho}(s) = R_{\mathrm{hol}}(s) - \sum_{j=1}^{J(s_0)} (-1)^{j}(s-s_0)^{-j}(\mathcal{L}_{X^{\rho}}+s_0)^{j-1}\Pi_k^{\rho}(s_0),
\end{equation}
where $R_{\mathrm{hol}}(s)$ is holomorphic near $s_0$, and $\Pi_k^{\rho}(s_0)$ is defined in \eqref{eq:Pi_k}. In addition,
$$\mathrm{WF}'(R_{\mathrm{hol}}(s))\subset \Delta(T^*M)\cup \Omega_+\cup E_u^*\times E_s^*, \quad \mathrm{WF}'(\Pi_k^{\rho}(s_0))\subset E_u^*\times E_s^*,$$
where $\Delta(T^*M)$ denotes the diagonal of $T^*M\times T^*M$, and $\Omega_+$ is the (positive) \emph{flow-out} defined by
$$\Omega_+:=\{(\varphi_t v,d\varphi_t^{-\top}\xi); (v,\xi)) ~:~ t\geq 0,\ \xi(X)=0\}. $$
Observe that $\WF([S_x\Sigma]\otimes  \pi^*\xi_\star) \subset \mathbb{H}^* \oplus \RR \alpha$. By \cite[Theorem 8.2.13]{Hor}, we have that $R_1^\rho(s)\big([S_x\Sigma]\otimes  \pi^*\xi_\star\big)$ is well-defined, meromorphic with respect to $s$, and
\begin{equation}
    \label{eq:WFresol}
    \mathrm{WF}\big(\Pi_1^{\rho}(0)\big([S_x\Sigma]\otimes  \pi^* \xi_\star\big)\big),\,\mathrm{WF}\big(R_{\mathrm{hol}}(s)\big([S_x\Sigma]\otimes  \pi^*\xi_\star\big)\big)\subset \Gamma,
\end{equation}
where
\[
\Gamma :=(\mathbb H^*\oplus \mathbb R \alpha)_x\cup  \bigcup_{t\geq 0}d\varphi_t^{-\top}\mathbb H_x^* \cup E_u^*, \qquad \mathbb{H}_x^* := \bigcup_{v \in S_x\Sigma} \mathbb H_v^*.
\]
In the above calculation, we have used that $(\mathbb H^*\oplus \mathbb R \alpha)_x\cap \{\xi(X)=0\}=\mathbb H^*_x$.

We claim that 
\begin{equation}
    \label{eq:WFdisjoint}
    \Gamma \cap d\varphi_{-\varepsilon}^{-\top}(\mathbb H^*\oplus \mathbb R \alpha)_y=\emptyset.
\end{equation}
First, $(\mathbb H^*\oplus \mathbb R \alpha)_x\cap d\varphi_{-\varepsilon}^{-\top}(\mathbb H^*\oplus \mathbb R \alpha)_y=\emptyset$ as $S_x\Sigma$ and $\varphi_{-\eps}(S_y\Sigma)$ are disjoint for small $\eps > 0$. Next, since $\mathbb{H}_x^* \subset\{\xi(X)=0\}$, the intersection $\big(\cup_{t\geq 0}(d\varphi_t^g)_x^{-\top}\mathbb H_x^*\big) \cap (d\varphi_{-\varepsilon}^g)_y^{-\top}(\mathbb H^*\oplus \mathbb R \alpha)_y$ is non-empty only when there is a $t\geq 0$ such that $(d\varphi_{t+\varepsilon}^g)_x^{-\top}\mathbb H_x^* \cap \mathbb H_y^*\neq \emptyset$ which is impossible by \eqref{eq:noconj}. Finally, $E_u^*$ is transverse to $(d\varphi_{-\varepsilon}^g)_y^{-\top}(\mathbb{H}^*\oplus \mathbb{R}\alpha)_y$. This proves \eqref{eq:WFdisjoint} and shows that the pairings \eqref{eq:pairing} and \eqref{eq:pairing2} have a meromorphic continuation to $\mathbb{C}$.
\end{proof} 

\subsection{Simplifying the pairing}
We now study the behavior of the Poincaré series at $s=0$.  

\begin{lem}
\label{lem:exactpt}
Let $y\in \Sigma$ and $\mu \in E_\rho^*(y)$ and extend $\mu$ to a smooth section $\mu_\star \in C^\infty(M;E_{\rho^*})$, parallel at $y$. Then $[S_y\Sigma]\otimes \pi^* \mu_\star$ is $d^{\nabla_{\rho^*}}$-exact. Moreover, the primitive can be chosen with wavefront set in $ (\mathbb H^*\oplus \mathbb R \alpha)_y.$
\end{lem}

\begin{proof} Since $\mu_\star$ is parallel at $y$, the current $[S_{y}\Sigma]\otimes \pi^*\mu_\star$ is $d^{\nabla_{\rho^*}}$-closed so it defines a twisted homology class in the corresponding twisted de Rham cohomology group $H^1_{\rho^*}(S\Sigma,\CC)$. By Poincaré duality, to show that this class is trivial, it suffices to show that:
\begin{equation*}
\forall \omega\in H^1_\rho(S\Sigma,\CC),\quad 
    \langle [S_{y}\Sigma]\otimes \pi^*\mu_\star, \omega \rangle= \int_{S_{y}\Sigma} \langle  \mu, \omega \rangle=0.
\end{equation*}
 By the Gysin exact sequence (see, for instance, the proof of \cite[Lemma 2.1]{HT26}), we may assume $\omega$ is the pullback of a form on the base $\Sigma$, i.e., $\omega =\pi^* \theta$ where $\theta \in H^1_\rho(\Sigma,\CC).$ Thus
$$\int_{S_{y}\Sigma} \langle  \mu, \omega \rangle=\int_{S_{y}\Sigma} \langle  \mu, \pi^*\theta \rangle=\int_{S_y\Sigma}\langle \mu,\pi^* i_y^*\theta\rangle=0,$$
where $i_y \colon \{y\}\to \Sigma$ is the inclusion of the point. That the primitive can be chosen with wavefront set contained in that of $[S_y\Sigma]\otimes \pi^* \mu_\star$ follows from the microlocal version of the Hodge theorem.
\end{proof}
From now on, $\xi_\star$ and $\mu_\star$ are sections on $\Sigma$ which coincide with $\xi$ and $\mu$ at $x$ and $y$ respectively, and are parallel at these points. We now simplify the expression in Proposition \ref{prop:DR} using the previous lemma.

\begin{prop}
\label{prop:primitive}
Let $(\Sigma,g)$ be a closed, connected, oriented, Anosov surface. Let $\rho\in \mathcal V_g \subset \mathcal{M}_r^{\mathrm{irr}}$, where $\mathcal V_g$ was defined in \eqref{def:Vg}. Let $F_{y,\mu}$ such that $d^{\nabla_{\rho^*}}F_{y,\mu} =[S_y\Sigma]\otimes \pi^*\mu_\star$. Then for any $\varepsilon>0$ small enough:
\begin{enumerate}
    \item If $x\neq y$, then
     \begin{equation}
     \label{eq:etaat0}
     \langle \eta_{g,\rho,x,y}(0)\xi,\mu\rangle =\langle [S_x\Sigma]\otimes \pi^* \xi_\star, e^{\varepsilon \mathcal L_{X^{\rho^*}}}F_{y,\mu} \rangle;
     \end{equation} 
      \item If $x= y$, then
     \begin{equation}
     \label{eq:etaat0bis}
     \langle \eta_{g,\rho,x,x}(0)\xi,\mu\rangle=\langle [S_x\Sigma]\otimes \pi^* \xi_\star, e^{\varepsilon \mathcal L_{X^{\rho^*}}}F_{x,\mu} \rangle+\langle \xi,\mu\rangle.
     \end{equation}
     \end{enumerate}
\end{prop}

\begin{rem}
    \label{rem:eps2}
  By Remark  \ref{rem:epsilon}, if $x\neq y$, we can take $\varepsilon=0$ in \eqref{eq:etaat0}.
\end{rem}
\begin{proof} It follows from Proposition \ref{prop:DR} that
     \begin{equation*}
     \begin{split}
\langle \eta_{g,\rho,x,y}(s)\xi,\mu\rangle-\delta_{x=y}\langle \xi,\mu\rangle &=e^{-s\varepsilon}\langle R_1^\rho(s)\big(\iota_X[S_x\Sigma]\otimes  \pi^*\xi_\star\big),e^{\varepsilon \mathcal L_{X^{\rho^*}}}\big([S_y\Sigma]\otimes  \pi^*\mu_\star\big)\rangle
\\&=e^{-s\varepsilon}\langle R_1^\rho(s)\big(\iota_X[S_x\Sigma]\otimes  \pi^*\xi_\star\big),e^{\varepsilon \mathcal L_{X^{\rho^*}}}d^{\nabla_\rho}F_{y,\mu}\rangle
\\&=e^{-s\varepsilon}\langle d^{\nabla_\rho}R_1^\rho(s)\big(\iota_X[S_x\Sigma]\otimes  \pi^*\xi_\star\big),e^{\varepsilon \mathcal L_{X^{\rho^*}}}F_{y,\mu}\rangle,
\end{split}
     \end{equation*}
where we used Lemma \ref{lem:exactpt} in the second line, $[d^{\nabla_\rho},e^{\varepsilon \mathcal L_{X^{\rho^*}}}]=0$ and an integration by parts in the third line. Continuing the computation, using $d^{\nabla_\rho}R_1^\rho(s)=R_2^\rho(s)d^{\nabla_\rho}$, we obtain:
\begin{equation}
\label{eq:1stpart}
\begin{split}
    \langle \eta_{g,\rho,x,y}(s)\xi,\mu\rangle-\delta_{x=y}\langle \xi,\mu\rangle  = & e^{-s\varepsilon}\langle R_2^\rho(s)d^{\nabla_\rho}\big(\iota_X[S_x\Sigma]\otimes  \pi^*\xi_\star\big),e^{\varepsilon \mathcal L_{X^{\rho^*}}}F_{y,\mu}\rangle
    \\= & e^{-s\varepsilon}\langle R_2^\rho(s)\mathcal L_{X^\rho}\big([S_x\Sigma]\otimes  \pi^*\xi_\star\big),e^{\varepsilon \mathcal L_{X^{\rho^*}}}F_{y,\mu}\rangle
    \\=&e^{-s\varepsilon}\langle \big([S_x\Sigma]\otimes  \pi^*\xi_\star\big),e^{\varepsilon \mathcal L_{X^{\rho^*}}}F_{y,\mu}\rangle \\
    & -se^{-s\varepsilon}\langle R_2^\rho(s)\big([S_x\Sigma]\otimes  \pi^*\xi_\star\big),e^{\varepsilon \mathcal L_{X^{\rho^*}}}F_{y,\mu}\rangle,
    \end{split}
\end{equation}
where we used Cartan's formula $\mathcal L_{X^\rho}=\iota_Xd^{\nabla_\rho}+d^{\nabla_\rho}\iota_X$ and Lemma \ref{lem:exactpt} which implies that $d^{\nabla_\rho}\big([S_x\Sigma]\otimes  \pi^*\xi_\star\big)=0$ in the second line, and in the last line the definition of the resolvent. The first term in the last line converges, when $s\to 0$, to the claimed formula. Hence, it is enough to show that the second term tends to $0$ as $s\to 0$.

By the same argument as in Proposition~\ref{prop:mero}, $\langle R_2^\rho(s)\big([S_x\Sigma]\otimes  \pi^*\xi_\star\big),e^{\varepsilon \mathcal L_{X^{\rho^*}}}F_{y,\mu}\rangle$ is meromorphic with respect to $s\in \mathbb{C}$. By \eqref{def:Vg}, for $\rho \in \mathcal V_g$, the space of resonant states on $2$-forms defined by \eqref{eq:res-1} has no Jordan blocks and is closed. Indeed, we have \[\Res^{2,1}(\rho,0)=\Res^{1,1}_0(\rho,0)\wedge \alpha,\] 
which follows from \eqref{eq:relalpha}. We can then write
\begin{equation*}
    R_2^{\rho}(s)=\Pi_2^{\rho}(0)/s+R_{\mathrm{hol}}(s) \colon C^{\infty}(M; \mathscr E_{\rho}^2)\to \mathcal{D}'(M; \mathscr E_{\rho}^2),\quad d^{\nabla_{\rho}} \Pi_2^{\rho}(0)=0,
\end{equation*}
where $\Pi_2^{\rho}(0)$ is defined in \eqref{eq:Pi_k} and $R_{\mathrm{hol}}(s)$ is holomorphic in $s$ in a neighborhood of $0$. Using this and the fact that $[S_x\Sigma]\otimes \pi^*\xi_\star= d^{\nabla_{\rho}}F_{x,\xi}$ is exact, we have
\[
\begin{split}
    \lim_{s\to 0} se^{-s\varepsilon}\langle R_s^{\rho}(s)([S_x\Sigma]\otimes \pi^*\xi_\star),e^{\varepsilon \mathcal{L}_{X^{\rho^*}}}F_{y,\mu}\rangle
    &= \lim_{s\to 0} e^{-s\varepsilon}\langle \Pi_2^{\rho}(0)([S_x\Sigma]\otimes \pi^*\xi_\star),e^{\varepsilon \mathcal{L}_{X^{\rho^*}}}F_{y,\mu}\rangle\\
    &=\langle \Pi_2^{\rho}(0) d^{\nabla_{\rho}}F_{x,\xi},e^{\varepsilon \mathcal{L}_{X^{\rho^*}}}F_{y,\mu}\rangle \\
    &=\langle d^{\nabla_{\rho}} \Pi_1^{\rho}(0) F_{x,\xi},e^{\varepsilon \mathcal{L}_{X^{\rho^*}}}F_{y,\mu}\rangle =0,
\end{split}
\]
where in the last line, we used the commutation relation $d^{\nabla_\rho}\Pi_1^\rho(0)=\Pi_2^\rho(0)d^{\nabla_\rho}.$
\end{proof}

\section{Computing the value at zero} 

\label{section:computation}

\subsection{Distinct points} In this section, we prove Theorem \ref{theo:1}. We start by the case of two distinct points.

\begin{prop}
    \label{prop:xneqy}
Let $(\Sigma,g)$ be a closed, connected, oriented, Anosov surface, and let $x,y \in \Sigma$ be distinct points. For any $\rho \in \mathcal V_g \subset \mathcal{M}_r$ $($where $\mathcal V_g$ is as in \eqref{def:Vg}$)$:
$$\eta_{g,\rho,x,y}(0)=0. $$
\end{prop}

We will use the following linear algebra lemma. Its proof is immediate and hinges on the fact that $\rho$ is irreducible if and only if the dual representation $\rho^*$ is irreducible.

\begin{lem} \label{lemma:linear-algebra}
    For $\rho \in \mathcal{M}_r^{\mathrm{irr}}$, $\mathrm{Span}\left\{\mathrm{Ran}(\rho(\gamma)-\mathrm{Id}) : \gamma \in \pi_1(\Sigma)\right\}=\mathbb C^r$.
\end{lem}

\begin{proof}[Proof of Proposition \ref{prop:xneqy}]
    We will construct an explicit primitive $F_{y,\mu}$ and apply Proposition \ref{prop:primitive}. For $\xi \in E_\rho(x)$, $\mu \in E_\rho^*(y)$, we will show that $\langle \eta_{g,\rho,x,y}(0)\xi,\mu\rangle=0$. It follows from Proposition \ref{prop:primitive} and Remark \ref{rem:eps2} that 
    $$\langle \eta_{g,\rho,x,y}(0)\xi,\mu\rangle=\langle [S_x\Sigma]\otimes \pi^* \xi_\star,F_{y,\mu} \rangle, $$
    where $\mathrm{WF}(F_{y,\mu})\subset (\mathbb H^*\oplus \mathbb R \alpha)_y. $ We will consider another primitive constructed as follows. 

    By Lemma~\ref{lemma:linear-algebra}, for $\rho\in \mathcal{V}_g$, there exists a finite number of elements $\gamma_j\in \pi_1(\Sigma)$ and $\theta_j \in E_\rho^*(y)$ such that $\sum_j(\rho^*(\gamma_j)\theta_j-\theta_j)=\mu$. Consider a parametrization $c_j \colon [0,1]\to \Sigma$ of a loop $c_j$ representing the class $\gamma_j$, based at $y$. Up to perturbing the curve, we can further assume that $x$ does not lie on $c_j$, for any $j$. For $t\in [0,1]$, let $\theta_{j,\star}(t)$ be the parallel transport of $\theta_j$ along $c_j$ from $y=c_j(0)$ to $c_j(t)$, and extend $\theta_{j,\star}$ locally transversally to $c_j$ by parallel transport. We then define
    \[
    \widetilde{F}_{y,\mu}:=\sum_j [S_{c_j}\Sigma]\otimes \pi^*\theta_{j,\star},
    \]
    where the above notation corresponds to the integration current over the restriction of $S\Sigma$ to $c_j$. Since $\theta_{j,\star}$ is parallel on the support of $S_{c_j}\Sigma$, 
    \begin{equation} \label{eq:primitive}
        d^{\nabla_{\rho^*}}\widetilde{F}_{y,\mu}=\sum_{j}[S_y\Sigma]\otimes \big((\rho(\gamma_j)-\mathrm{Id})\theta_j\big)=[S_y\Sigma]\otimes \pi^*\mu_\star.
    \end{equation}
Hence $\widetilde{F}_{y,\mu}$ is another primitive of $[S_y\Sigma]\otimes \pi^* \mu_\star$. Note that the support of $\widetilde{F}_{y,\mu}$ is included in $\cup_{j=1}S_{c_j}\Sigma$, the union of the unit tangent bundles restricted to $c_j$. In particular, its wavefront set is included in the conormal bundles $\cup_{j=1}N^*S_{c_j}\Sigma$. Since the wavefront set of $F_{y,\mu}$ is included in $N^*(S_y\Sigma)$, the following holds:
$$d^{\nabla_{\rho^*}}(F_{y,\mu}-\widetilde{F}_{y,\mu})=0,\quad \mathrm{WF}(F_{y,\mu}-\widetilde{F}_{y,\mu})\subset N^*S_y\Sigma \cup \bigcup_{j=1}N^*S_{c_j}\Sigma=:\mathcal C.  $$
By the Hodge decomposition theorem (see \cite[Lemma 2.1]{zazi}), there is a $d^{\nabla_{\rho^*}}$-closed 1-form $\omega \in \Omega^1 \otimes \mathcal E_{\rho^*}$ and $f\in \mathcal D'(M;\mathcal E_{\rho^*})$ with $\WF(f)\subset \mathcal{C}$ such that 
    $$F_{y,\mu}-\widetilde{F}_{y,\mu}=\omega+d^{\nabla_{\rho^*}}f. $$
Since $x$ does not lie on any of the curves $c_j$, $\mathcal C$ is disjoint from the wavefront set of $[S_x\Sigma]\otimes \pi^* \xi_\star$. 
    In particular, the pairing can be rewritten as
$$\langle \eta_{g,\rho,x,y}(0)\xi,\mu\rangle = \langle [S_x\Sigma]\otimes \pi^* \xi_\star,\widetilde{F}_{y,\mu} \rangle+\langle [S_x\Sigma]\otimes \pi^* \xi_{\star},\omega \rangle+\langle [S_x\Sigma]\otimes \pi^*  \xi_\star,d^{\nabla_{\rho^*}}f \rangle.$$
Since $[S_x\Sigma]\otimes \pi^* \xi_\star$ is supported above $x$ and $\widetilde{F}_{y,\mu}$ is supported on the curves $c_j$, the first pairing on the right-hand side is zero. Regarding the second pairing, we use Lemma \ref{lem:exactpt} and the fact that $\omega$ is closed to obtain
$$\langle [S_x\Sigma]\otimes \pi^* \xi_\star,\omega \rangle=\langle d^{\nabla_{\rho}}F_{x,\xi},\omega \rangle=\langle F_{x,\xi} ,d^{\nabla_{\rho^*}}\omega \rangle=0. $$
To deal with the third pairing, we use Lemma \ref{lem:exactpt} and integrate by parts to obtain
$$\langle [S_x\Sigma]\otimes \pi^* \xi_\star,d^{\nabla_{\rho^*}}f \rangle=  \langle d^{\nabla_\rho}([S_x\Sigma]\otimes \pi^* \xi_\star),f \rangle=0. $$
This proves that the Poincaré series vanishes at $s=0$.
\end{proof}

\subsection{Equal points} We now treat the case of two equal points.

\begin{prop}
    \label{prop:xeqy}
Let $(\Sigma,g)$ be a closed, connected, oriented, Anosov surface, and let $x \in \Sigma$. For any $\rho \in \mathcal V_g \subset \mathcal{M}_r$, where $\mathcal V_g$ is the subset of representations from \eqref{def:Vg}, we have:
$$\eta_{g,\rho,x,x}(0)=0. $$
\end{prop}

\begin{proof}
    As in the proof of Proposition \ref{prop:xneqy}, we will work with a more explicit primitive. Note that $\pi(\varphi_{-\varepsilon}(S_x\Sigma))$ bounds a simply connected open set $D_\varepsilon(x)$ containing $x$ (the geodesic disc of radius $\varepsilon>0$ centered at $x$). Since $E_\rho$ is flat, we can further assume that $\mu_\star$ (such that $\mu_\star(x)=\mu$) is parallel on this disc. To construct an explicit primitive of $e^{\varepsilon \mathcal L_{X^\rho}}([S_x\Sigma]\otimes \pi^*\mu_\star)$, we first choose a point $z\in D_\varepsilon(x)\setminus\{x\}$ such that $d(x,z)<\varepsilon/3$. Next, we define a vector field $Y$ on $D_{\varepsilon}(x)\setminus \{z\}$ as follows:
    \begin{enumerate}
        \item If $w\in D_{\varepsilon/2}(z)\setminus \{z\}$, write $w=\pi(\varphi_{-s} u)$ for $s\in [0,\varepsilon/2)$, $u\in S_z\Sigma$, and let $Y(\pi(\varphi_{-s} u)):= \varphi_{-s} u$.
        \item If $w\in \pi(\varphi_{-\varepsilon}^g(S_x\Sigma))$, write $w=\pi(\varphi_{-\varepsilon} v)$ for $v\in S_x\Sigma$, and let $Y(\pi(\varphi_{-\varepsilon} v)):= \varphi_{-\varepsilon} v$.
        \item The vector field $Y$ is smooth on $D_{\varepsilon}(x)\setminus \{z\}$ with a singularity at $z$.
    \end{enumerate}
    The two prescribed vector fields on the boundary have degree $1$ so they are homotopic through a non-vanishing vector field on the annulus and thus such a vector field $Y$ exists. We then define a surface with boundary in $S\Sigma$:
    \begin{equation}
        \label{eq:fancyY}
        \mathcal Y:=\{(w,Y(w)) ~:~ w\in D(x,\varepsilon)\setminus \{z\}\}.
    \end{equation}
    We then consider $[\mathcal Y]\otimes \pi^*\mu_\star$. Since $\pi^* \mu_\star$ is parallel with respect to the pullback flat connection on $\mathcal{E}_\rho = \pi^*E_\rho$, it holds that:
    \[
    d^{\nabla_\rho}([\mathcal Y]\otimes \pi^*\mu_\star)=e^{\varepsilon \mathcal L_{X^\rho}}([S_x\Sigma]\otimes \pi^*\mu_\star)-[S_z\Sigma]\otimes  \pi^*\mu_\star.
    \]
    Here, we have oriented $\mathcal Y$ by the orientation induced by the inward normal on $\pi(\varphi_{-\varepsilon}(S_x\Sigma))$.
    We can write
    \begin{equation}
    \label{eq:goodprimitive}e^{\varepsilon \mathcal L_{X^\rho}}([S_x\Sigma]\otimes \pi^*\mu_\star)=d^{\nabla_\rho}\big([\mathcal Y]\otimes \pi^*\mu_\star+\widetilde{F}_{z,\mu} \big), \end{equation}
    where $\widetilde{F}_{z,\mu}$ is the primitive constructed in the proof of Proposition \ref{prop:xneqy}.

In order to use this explicit primitive in \eqref{eq:etaat0}, we now compute the wavefront set of $[\mathcal Y]\otimes \pi^*\mu_\star$. Since we will pair it with the current $[S_x\Sigma]\otimes \pi^*\mu_\star$, which is supported at $x$, we only need to compute the wavefront set of $[\mathcal Y]\otimes \pi^*\mu_\star$ near $x$. Given our choice of $z$, we see that near $x$, $Y$ is defined by item (1). Thus, $\mathcal Y$ is parametrized near $x$ by $\varphi_{-s}(u)$ for some $s\in (0,\varepsilon/2)$ and $u\in S_z\Sigma.$ We compute the tangent space to $\mathcal Y$ at the point $(x,Y(x))=\varphi_{-s_0}(u_0)$. Observe that 
    $$\partial_s\varphi_{-s}(u_0)|_{s=s_0}= - X(x, Y(x)), \quad \partial_u\varphi_{-s_0}(u)|_{u=u_0}=d(\varphi_{-s_0})_{(z,u_0)}(V). $$
    Therefore
    $$T_{(x,Y(x))}\mathcal Y=\mathbb R X(x,Y(x)) \oplus d(\varphi_{-s_0})_{(z,u_0)}(\mathbb V). $$
    Note that by the choice of orientation of $\mathcal Y$, we deduce that a positively oriented basis of $T_{(x,Y(x))}\mathcal Y$ is given by $(-X(x,Y(x)),d(\varphi_{-s_0})_{(z,u_0)}(V(z,u_0)).$ 
    Using \eqref{eq:positive}, this implies
    $$T_{(x,Y(x))} S_x\Sigma \oplus T_{(x,Y(x))}\mathcal Y  =T_{(x,Y(x))} M, $$
    so the conormal bundles to $\mathcal Y$ and $S_x\Sigma$ are  transverse at their unique point of intersection. In particular, $[\mathcal Y]\otimes \pi^*\mu_\star+ \widetilde{F}_{z,\mu}$ is a primitive of $e^{\varepsilon \mathcal L_{X^\rho}}([S_x\Sigma]\otimes \pi^*\mu_\star)$ with wavefront set disjoint from that of $[S_x\Sigma]\otimes \pi^*\mu_\star$, and we deduce from an argument analogue to the one used in the proof of Proposition \ref{prop:xneqy} that we can use this primitive in \eqref{eq:etaat0bis}. This yields
\begin{align*}
   \langle \eta_{g,\rho,x,x}(0)\xi,\mu\rangle&= \langle [S_x\Sigma]\otimes \pi^*\xi_\star ,[\mathcal Y]\otimes \pi^*\mu_\star\rangle +\langle[S_x\Sigma]\otimes \pi^*\xi_\star,\widetilde{F}_{z,\mu}\rangle+\langle \xi,\mu\rangle
   \\&=\langle [S_x\Sigma]\otimes \pi^*\xi_\star ,[\mathcal Y]\otimes \pi^*\mu_\star\rangle+\langle \xi,\mu\rangle
   \\&=\langle [S_x\Sigma], [\mathcal Y]\rangle \times \langle \xi,\mu\rangle+\langle \xi,\mu\rangle,
\end{align*}
where we used in the second equality that $\langle[S_x\Sigma]\otimes \pi^*\xi_\star,\widetilde{F}_{z,\mu}\rangle$ vanishes by the previous proof since $z\neq x$, and in the last equality that $\xi_\star$ and  $\mu_\star$ are parallel. The proof will be complete after we prove $\langle [S_x\Sigma], [\mathcal Y]\rangle = -1$. To do so, let us remark that $S_x\Sigma$
and $\mathcal Y$ intersect at a single point. The pairing is thus equal to $\pm 1$, with the sign being determined by the orientations.

It follows from the previous discussion that an oriented basis of $T_{(x,Y(x))} S_x\Sigma \oplus T_{(x,Y(x))}\mathcal Y$ is given by 
$$\big(V(x,Y(x)), -X(x,Y(x)),d(\varphi_{-s_0})_{(z,u_0)}(V(z,u_0))\big). $$
Applying \eqref{eq:positive}, we see that this basis has the same orientation as
$$\big(V(x,Y(x)),-X(x,Y(x)),a_{-s_0}(z,u_0)H(x,Y(x)))\big),$$
with $a_{-s_0}<0$. From \eqref{eq:orientation}, we see that this basis is negatively oriented. This implies that $\langle [S_x\Sigma], [\mathcal Y]\rangle=-1$, which concludes the proof.
\end{proof}

\subsection{Alternative proof}
    We sketch here an alternative proof of Theorem \ref{theo:1}. We claim that for any $x\in \Sigma$, any $\xi \in (E_\rho)_x$ and any $\rho \in \mathcal V_g$,
$$y\mapsto \eta_{g,\rho,x,y}(0)\xi  \text{ is parallel, i.e., } \nabla_y \eta_{g,\rho,x,y}(0)\xi=0.  $$ Indeed, for $\Re(s) \gg 1$, the function $y\mapsto \eta_{g,\rho,x,y}(s)\xi$ is $C^1$ except at $y=x$. Differentiating \eqref{eq:Poinc}, gives
   $$\forall x\neq y, \quad \nabla_y \eta_{g,\rho,x,y}(s) = - s \sum_{\gamma \in \mathcal{P}_{x,y}} \nabla_y \ell_g(\gamma) e^{-s\ell_g(\gamma)}\tau_\rho(\gamma), $$
   where we used that $\tau_\rho(\gamma)$ is defined by parallel transport. Note that we have $\nabla_y \ell_g(\gamma) = \dot{\gamma}(\ell_g(\gamma))$. Let $\beta \in C^\infty(T^*M)$. Using the same arguments as in \eqref{eq:long}, we can then write, for $\xi \in (E_\rho)_x$ and $\mu \in (E_{\rho^*})_y$, 
\[
\begin{split}
\langle\nabla_y \eta_{g,\rho,x,y}(s)\xi, \beta(y)\otimes \mu\rangle& = - s \sum_{\gamma \in \mathcal{P}_{x,y}} \langle \nabla_y \ell_g(\gamma) e^{-s\ell_g(\gamma)}\tau_{\rho} (\gamma) \xi, \beta(y) \otimes \mu\rangle \\
& = - s \langle R_{1,\rho}^{\pi^*T\Sigma}(s)(\iota_X[S_x\Sigma]\otimes \pi^* \xi_\star \otimes \mathbf{s}), [S_y\Sigma]\otimes \pi^* \mu_\star \otimes \pi^*\beta\rangle,
\end{split}
\]
where $\mathbf{s} \in C^\infty(S\Sigma, \pi^*T\Sigma)$, $\mathbf{s}(x,v) := v$ is the tautological section and $R_{1,\rho}^{\pi^*T\Sigma}(s)$ is the resolvent on $\mathscr E_{\rho,0}^1 \otimes \pi^*T\Sigma$, where $X$ acts by $\mathcal{L}_{X^\rho}\otimes(\pi^*\nabla^{\mathrm{LC}})_X$. Observe that $(\pi^*\nabla^{\mathrm{LC}})_X\mathbf{s} = 0$ so the resolvent $R_{1,\rho}^{\pi^*T\Sigma}(s)$ acts on a pure state $\omega \otimes s$ $($where $\omega \in \mathcal{D}'(M; \mathscr E_{\rho,0}^1))$ as $R_{1,\rho}^{\pi^*T\Sigma}(s) (\omega \otimes \mathbf{s}) = (R_{1}^\rho(s)\omega) \otimes \mathbf{s}$, where $R_{1}^\rho(s)$ is the resolvent on $1$-forms. We thus find:
\begin{equation}
\label{eq:equal}
\langle\nabla_y \eta_{g,\rho,x,y}(s)\xi , \beta(y)\otimes \mu\rangle = -s \langle R_1^\rho(s)(\iota_X[S_x\Sigma]\otimes \pi^*  \xi_\star) \otimes \mathbf{s}, [S_y\Sigma]\otimes \pi^*  \mu_\star \otimes \pi^*\beta \rangle.
\end{equation}
From Proposition \ref{prop:HT} and Lemma \ref{lem:exactpt}, we have
\begin{align*}R_1^\rho(s)(\iota_X[S_x\Sigma]\otimes \pi^*  \xi_\star)&=R_1^{\rho,\mathrm{hol}}(s)(\iota_X[S_x\Sigma]\otimes \pi^*  \xi_\star)+s^{-1}\Pi_1^\rho(0)(\iota_X d^{\nabla_\rho}F_{x,\xi})
\\&=s^{-1}\underbrace{\Pi_1^\rho(0)\mathcal L_{X^\rho}}_{=0} F_{x,\xi}- s^{-1}\underbrace{d^{\nabla_\rho}\Pi_0^\rho(0)}_{=0}\iota_X F_{x,\xi}+O(1)=O(1).
\end{align*}
In particular, \eqref{eq:equal} implies that $\nabla_y \eta_{g,\rho,x,y}(0)=0$ when $y\neq x$. Our choice to include the trivial loop in the Poincaré series when $x=y$ makes $y\mapsto \eta_{g,\rho,x,y}(0)$ continuous at $y=x$ and hence this proves the claim for any $x,y\in \Sigma$. Now, let $p\in \Sigma$ be a periodic point, then 
$$\forall \gamma \in \pi_1(\Sigma),\quad \rho(\gamma) \eta_{g,x,p,\rho}(0)\xi=\eta_{g,x,p,\rho}(0) \xi. $$
Using Lemma \ref{lemma:linear-algebra}, we deduce that that $\eta_{g,x,p,\rho}(0)\xi=0$. This proves that the whole section vanishes by the density of periodic points and the fact that the section is continuous.

\section{Proof of Theorems \ref{theo:2} and \ref{theo:3}}

\label{sec:7}

\subsection{Resonant forms at zero} We start by describing the structure of resonant states at zero for a hyperbolic surface $(\Sigma,g)$ and the representation associated to a holomorphic $1$-form $\omega$ close enough to $0$ (Theorem \ref{theo:3}). 

\subsubsection{Kernel of the twisted Laplacian} The \emph{twisted flat Laplacian} on $C^\infty(\Sigma;E_{\rho_\omega})$ is
\[
\Delta_{\omega}:=-\operatorname{tr}\big((\nabla^{\rho_\omega})^2\big).
\]
Here, the second covariant derivative is the induced derivation $\nabla^{\rho_\omega} \colon C^\infty(\Sigma;T^*\Sigma \otimes E_{\rho_\omega}) \to C^\infty(\Sigma;T^*\Sigma^{\otimes 2} \otimes E_{\rho_\omega})$, where the Levi-Civita connection is applied to the $T^*\Sigma$ factor. Let $p \colon \mathbb{H}^2 \to \Sigma = \mathbb{H}^2/\Gamma$ be the footpoint projection (of course, this $\Gamma$ is different from the one in \eqref{eq:WFresol}). In the trivialization of
\(p^*E_{\rho_\omega} = \CC \times \mathbb{H}^2\) over \(\mathbb H^2\), $p^*C^\infty(\Sigma;E_{\rho_\omega})$ identifies with the space of equivariant sections $f \in C^\infty(\mathbb{H}^2,\CC)$ such that
\begin{equation}
    \label{equation:equivariance}
f(\gamma.x) = \exp\left(\int_\gamma \omega \right) f(x).
\end{equation}
The lift of $\Delta_\omega$ to $\mathbb{H}^2$ is the scalar hyperbolic Laplacian. In general, the twisted Laplacian is \emph{not} selfadjoint. 

\begin{lem}
 \label{lemm:gap}
     There is a neighborhood $U$ of $0$ in $H^0(\Sigma;K)$ such that for any $\omega \in U$, $0$ is an eigenvalue of $\Delta_{\omega}$ of algebraic  multiplicity $1$.
 \end{lem}

\begin{proof}
Let $F$ be a holomorphic primitive of $p^*\omega$ on the universal cover $\mathbb H^2$. $e^F$ satisfies \eqref{equation:equivariance}.
Since $\omega$ is holomorphic, $e^F$ is harmonic: $\Delta_{\mathbb H^2} e^F = 0$. This proves that $\mathrm{dim}(\ker \Delta_\omega)\geq 1$ by the previous discussion.  Since $\ker \Delta_0$ is $1$-dimensional, the same holds for all $\omega$ close to~$0$.
\end{proof}

\subsubsection{Resonant $0$- and $2$-forms} We describe them using Lemma \ref{lemm:gap} and the classical-quantum correspondence of \cite{GHW}.

\begin{lem}
    \label{lem:CQ}
    There is a neighborhood $U$ of $0$ in $H^0(\Sigma;K)$ such that for any $\omega \in U,$ we have
    \begin{equation}
        \label{eq:Res02}
        \mathrm{dim}(\Res_0^{k,\infty}(\rho_{\omega},0))=1, \quad k=0,2.
    \end{equation}
\end{lem}

\begin{proof}
    By \cite[Proposition 3.1]{zazi}, \eqref{eq:Res02} holds at $\omega=0$. By perturbation theory \cite{B20}, if $U \subset H^0(\Sigma;K)$ is chosen small enough around $0$, then for any $\omega \in U,$ 
    $$\mathrm{dim}(\Res_0^{k,\infty}(\rho_{\omega},0))\leq 1, \quad k=0,2. $$
    Since $\mathrm{dim}(\Res_0^{0,\infty}(\rho_{\omega},0))=\mathrm{dim}(\Res_0^{2,\infty}(\rho_{\omega},0)),$ we only need to show that $\Res_0^{0,1}(\rho_{\omega},0)$ is non-trivial. The idea is to show that one can construct from the eigenvector $e^F$ (where $F$ is a holomorphic primitive of $p^*\omega$ on the universal cover $\mathbb H^2$) of $\Delta_\omega$ associated to $0$, given by Lemma \ref{lemm:gap}, a resonant state using the classical-quantum correspondence.
    
     Note that for any $o \in \mathbb{H}^2$, there exists $C >0$ such that $|e^{F(x)}| \leq e^{Cd(x,o)}, \forall x \in \mathbb{H}^2$. We can thus apply the \emph{Poisson--Helgason transform} $\mathcal P_0$ to $e^F$. The exact definition of $\mathcal P_0$ will not be needed and can be found in \cite[Lemma 2.1]{GHW}. From \cite[Lemma 2.1]{GHW}, we obtain that for the value $\lambda=0$, the Poisson--Helgason transform is an isomorphism from $\mathcal D'(\mathbb S^1)$ to tempered distributions in the kernel of $\Delta_{\mathbb H^{2}}$. Let $\Xi \in \mathcal D'(\mathbb S^1)$ be the unique distribution such that $\mathcal P_0 (\Xi)=e^F$. Let $B_- \colon S\mathbb H^2\to \mathbb S^1$ be the backward endpoint map in the Hopf coordinates, and let 
    \begin{equation}
        \label{eq:tildeT_F}
       \widetilde \Xi(x,v):=B_-^*\Xi(x,v),\quad  \widetilde{\Xi}_{F}(x,v):=e^{-F(x)} \widetilde \Xi(x,v), \qquad (x,v) \in S\mathbb{H}^2
    \end{equation}
    where the pullback is well defined since $B_-$ is a submersion. By \cite[Equation $(2.12)$]{GHW}: 
    \begin{equation}\label{eq:td_Xi}
        X\widetilde{\Xi}=0, \quad U_-\widetilde \Xi=0, \quad \mathrm{WF}(\widetilde \Xi)\subset  E_u^*, 
    \end{equation}
    where $U_-$ denotes the unstable horospherical vector field \eqref{equation:horocyclic}. See the proof of \cite[Theorem 3.3]{GHW}. Since $F$ is a smooth primitive of $p^*\omega$, we find that
    $$(X+\omega(dp(v)))\widetilde \Xi_F(x,v)=0,\quad \mathrm{WF}(\widetilde \Xi_F)\subset  E_u^*.$$ Finally, by \cite[Lemma 2.1]{GHW}, $\widetilde \Xi$ is equivariant under the action of $\Gamma$ from which we deduce that $\widetilde \Xi_F$ satisfies \eqref{equation:equivariance}. This shows that it descends to a (nonzero) distribution on $S\Sigma$ with values in the flat bundle and, by the previous discussion, it defines a nontrivial twisted resonant state at zero.
\end{proof}

We write $\mathrm{Res}^{0,\infty}_0(\rho_{\omega}, 0)=\mathrm{Span}(h)$, and record further information about $h$. 

 \begin{lem}
 \label{lem:renormalization}
     Write the Fourier decomposition $h=\sum_{k\in \mathbb Z}h_k$ with $h_k\in \Omega_k$ $($see \S\ref{sssection:fourier}$)$. Then we can normalize $h$ so that 
     $$h_{-1}=0,\quad h_0=1,\quad h_1=-2\omega.$$
 \end{lem}

 We will use the operators introduced in \S\ref{sssection:twisted-geometry} for the proof.

 \begin{proof}
     The twisted Casimir operator $C^\omega$ (see \eqref{eq:Casimir}) commutes with $V$ (Lemma \ref{lemma:casimir}) and satisfies $C^\omega|_{\Omega_0}=-\Delta_\omega$. 
     Since $h\in \mathrm{Res}^{0,\infty}_0(\rho_{\omega}, 0)$ and $\mathrm{Res}^0(\rho_\omega)\subset \{z\in \CC, \mathrm{Re}(z)\leq \|\omega\|_{{\infty}}\}$, \eqref{eq:td_Xi} shows that $h$ is a first band resonant state (i.e., it lies in the kernel of $U_{-}^{\omega}$) provided $\|\omega\|_{{\infty}}$ is small enough. Therefore, $X^\omega h=0= U_-^\omega h$. In particular, by \eqref{eq:Casimir}, $-\Delta_\omega h_0=C^\omega h_0=(C^\omega h)_0=0$.
This means that $h_0\in \mathrm{Ker}(\Delta_\omega)$ and by Lemma \ref{lemm:gap}, we can suppose that $h_0=1$ (we have $h_0\neq 0$ since by \eqref{eq:eta_omega} below, $h_0=0$ implies that $h_k=0$ for all $k\in \mathbb{Z}$).

An adaptation of \cite[Lemma 2.2]{GHW} shows that for any $k\in \mathbb Z$, we have (with $\eta^\omega_- = \eta_-$):
     \begin{equation}\label{eq:eta_omega}
         2\eta^\omega_\pm h_{k\mp 1}=\mp k h_k. 
     \end{equation} With $k=1,$ this gives $h_{1}=-2\eta^\omega_+(1)=-2\omega$ and with $k=-1$, we obtain $h_{-1}=-2\eta_-(1)=0.$ This concludes the proof.
  \end{proof}

\subsubsection{Resonant $1$-forms}

We now study the generalized resonant states for the action on one $1$-forms.

\begin{prop} \label{prop:gen-res-states}
Let $(\Sigma,g)$ be a hyperbolic surface of genus $\textsl{g}\geq 2$. There exists a neighborhood $U$ of $0$ in $H^0(\Sigma; K)$, such that for any $\omega \in U\setminus\{0\}$, we have
\begin{equation}
    \label{eq:Res1}
    \mathrm{dim}(\Res_0^{1,\infty}(\rho_{\omega},0))=2\textsl{g},\quad \mathrm{dim}(\Res_0^{1,1}(\rho_{\omega},0))=2\textsl{g}-1.
\end{equation}
Moreover, we have $d^{\nabla_{\rho_\omega}}(\Res_0^{1,1}(\rho_{\omega},0))=0$. That is, all resonant states are $d^{\nabla_{\rho_\omega}}$-closed.
\end{prop}

\begin{proof}
    By \cite[Corollary C]{FS}, the order of vanishing of the zeta function of $(g,\rho_\omega)$ at zero is equal to $2\textsl{g}-2$. By Lemma \ref{lem:CQ} and \cite[Equation (1.2)]{HT26} we deduce that we have $\mathrm{dim}(\Res_0^{1,\infty}(\rho_{\omega},0))=2\textsl{g}.$ We now study more precisely the Jordan block structure at zero. Since $g$ is hyperbolic, the vector bundle of $1$-forms in $\ker \iota_X$ is trivialized over $S\Sigma$ by two smooth sections $\beta_\pm$ such that $\mathcal{L}_X\beta_\pm =\pm \beta_\pm$. These $1$-forms are dual to the stable and unstable horocyclic vector fields. Moreover, we have
 $$\mathcal L_{X^{\rho_\omega}}(f_+\beta_++f_-\beta_-)=(X^\omega(f_+)+1)\beta_++(X^\omega(f_-)-1)\beta_-. $$
 This means that
 $$\mathrm{Res}_0^{1,\infty}(\rho_{\omega},0)\cong \mathrm{Res}_0^{0,\infty}(\rho_{\omega},1)\oplus \mathrm{Res}_0^{0,\infty}(\rho_{\omega},-1). $$
 Note that the isomorphism preserves the Jordan block structure. Since $\mathrm{Res}^0(\rho_\omega)\subset \{z\in \CC, \mathrm{Re}(z)\leq \|\omega\|_{{\infty}}\},$ if $\|\omega\|_\infty$ is small enough, then $\mathrm{Res}_0^{0,\infty}(\rho_{\omega},1)=0$ and it remains to study $\mathrm{Res}_0^{0,\infty}(\rho_{\omega},-1)$.
 Using Lemma \ref{lemma:casimir}, we have
 \begin{equation}
 \label{eq:Umap}U_-^\omega \colon \mathrm{Res}^{0,\infty}_0(\rho_{\omega},-1)\to \mathrm{Res}^{0,\infty}_0(\rho_{\omega}, 0), 
 \end{equation}
 where the space on the right is $1$-dimensional by Lemma \ref{lem:CQ}.

  \emph{Step 1. Every resonant state is in the first band.} Given $f\in \mathrm{Res}^{0,1}_0(\rho_{\omega},-1)$, we now show that $f\in \mathrm{Ker}(U_-^\omega).$ By \eqref{eq:Umap}, there is a $c\in \mathbb C$ such that 
  $$(X^\omega-1)f=0,\quad  U_-^\omega f=ch,$$
  where $h$ is normalized as in Lemma \ref{lem:renormalization}. Projecting onto $\Omega_k$, this implies that for any $k\in \ZZ$
   \begin{equation} \label{eq:nu-pm-modes}
       2\eta_{\pm}^\omega f_{k\mp 1}=(1\mp k)f_k\mp ich_k.
   \end{equation}
   With $k=-1,$ we obtain $2\eta_-f_0=ich_{-1}=0$ so $f_0$ is constant. Applying \eqref{eq:nu-pm-modes} with $k=1$ gives
   $$2\omega f_0 =2\eta_+^\omega(f_0)=-ic h_1=2ic\omega,$$
   where we used Lemma \ref{lem:renormalization} again. Since $\omega\neq 0,$ this implies $f_0=ic$. Integrating \eqref{eq:nu-pm-modes} with $k=0$ and $\eta_-$ gives, using Lemma \ref{lem:renormalization} and integration by parts, that $f_0$ integrates to zero and hence $c=0,$ that is,  $f\in \mathrm{Ker}(U_-^\omega).$
\medskip

\emph{Step 2. Every resonant state is $d^{\nabla_{\rho_\omega}}$-closed.} The lifted twisted exterior derivative on $S\Sigma$ is given by $d+\pi^*\omega \wedge$. Observe that $\pi^*\omega(X)=\omega_x(v)$ and $\pi^*\omega(H)=\omega_x(Jv) = i \omega_x(v)$ since $\omega$ is of type $(1,0)$. The hyperbolic structural equations give $d\beta_-=-\alpha \wedge \beta_-.$ Given $f\in \mathrm{Res}^{0,1}_0(\rho_{\omega},-1)$, we then find:
\begin{align*}d^{\nabla_{\rho_\omega}}(f\beta_-)&=d^{\nabla_{\rho_\omega}}f\wedge \beta_-+ fd \beta_-=(X^\omega(f)\alpha+ U_-^\omega(f) \beta_+)\wedge \beta_--f\alpha \wedge \beta_-
\\&=\big( (X^\omega-1)f\big)\alpha\wedge \beta_-+(U^\omega_- f) \beta_+\wedge \beta_-=0, 
\end{align*}
by the first step. This shows that all resonant states are closed.

\emph{Step 3. The space of resonant states is at most $(2\textsl{g}-1)$-dimensional.} 
By Step $2$, the elements of $\Res_0^{1,1}(\rho_{\omega},0)$ are closed. This implies that there is a natural map (see \cite[Lemma 4.3]{HT26} for the exact definition):
\begin{equation}\label{eq:Res0beta_map}
   S : \Res_0^{1,1}(\rho_{\omega},0) \to H^1(S\Sigma,\rho_{\omega}),
\end{equation}
where $H^1(S\Sigma,\rho_{\omega})$ is the twisted cohomology group of $S\Sigma$. 
By \cite[Lemma 2.1]{HT26}, $\dim H^1(S\Sigma,\rho_{\omega}) =2\textsl{g}-2$ and the kernel of $S$ is equal to $d^{\nabla_{\rho_{\omega}}}(\Res_0^{0,1}(\rho_{\omega},0))=\mathbb{C} d^{\nabla_{\rho_{\omega}}}h$. Therefore, applying the rank--nullity theorem to $S$ yields
\begin{equation} \label{eq:rn-res}
   \mathrm{dim}( \Res_0^{1,1}(\rho_{\omega},0))\leq \mathrm{dim}(H^1(S\Sigma,\rho_{\omega}))+\mathrm{dim}(\mathrm{Ker}(S))=2\textsl{g}-2+1=2\textsl{g}-1.
\end{equation}

\emph{Step 4. The space of resonant states is $(2\textsl{g}-1)$-dimensional.} We claim that $h d\alpha$ is exact. First,  by the twisted Gysin exact sequence in \cite[Lemma 2.1]{HT26}, we have the isomorphism
 \begin{equation*}
   \pi_* \colon  H^2(S\Sigma,\rho_{\omega}) \to H^1(\Sigma,\rho_{\omega}).
 \end{equation*}
We now compute $\pi_*(h d\alpha)$. For any smooth $1$-form $\theta$ on $\Sigma,$ a direct computation using the structural equations gives $d\alpha \wedge \pi^*\theta=\theta_x(v)\alpha \wedge d\alpha$. This means that 
\begin{align*}\int_{\Sigma}\pi_*(h d\alpha)\wedge \theta =\langle \pi_*(h d\alpha),\theta\rangle=\langle h d\alpha,\pi^*\theta\rangle=\int_{S\Sigma} h(x,v) \theta_x(v) \alpha \wedge d\alpha= - \pi \int_{\Sigma} \star h_1\wedge \theta ,  
\end{align*}
where $\star$ is the Hodge star on $1$-forms.
We hence deduce that $\pi_*(h d\alpha)=\star h_1.$ Lemma \ref{lem:renormalization} and the fact that $\star \omega=-i\omega$ imply that $\pi_*(h d\alpha)$ is a multiple of $\omega$. It is hence exact in $H^1(\Sigma,\rho_{\omega})$ as $d^{\nabla_{\rho_\omega}}\cdot 1=\omega$.

 The same argument applied to the dual representation $\rho_{-\omega}$ gives the exactness of the co-resonant state $h_*d\alpha\in \Res_0^{*,2,1}(\rho_{-\omega},0)$. By the argument of \cite[Lemma 2.7]{CDDP}, the map~\eqref{eq:Res0beta_map} is surjective. Thus, the inequality in \eqref{eq:rn-res} is an equality, and we conclude
 \begin{equation*}
      \dim(\Res_0^{1,1}(\rho_{\omega},0) )  = 2\textsl{g}-1.\qedhere
 \end{equation*}
 \end{proof}

\begin{proof}[Proof of Theorem \ref{theo:3}] Follows from Lemma \ref{lem:CQ} and Proposition \ref{prop:gen-res-states}.
\end{proof}

\subsection{Value at zero}

We now use the structure of resonant states obtained in Proposition \ref{prop:gen-res-states} to show that the Poincaré series admits a well-defined value at zero, and we compute it. In the rest of the proof, we will use the subscript $\omega$ instead of $\rho_\omega$ to denote objects corresponding to the $1$-dimensional representation associated to $\omega$.

\subsubsection{Complex geometry} Let $H^0(\Sigma;E_\omega)$ be the space of holomorphic sections of $E_\omega \to \Sigma$.

\begin{lem}
\label{lemma:evaluation}
    For $\omega \in H^0(\Sigma;K)$ small enough, $H^0(\Sigma;E_\omega) \simeq \CC$ and for all $x \in \Sigma$, the evaluation map
    \[
    \mathrm{ev}_x^\omega \colon H^0(\Sigma;E_\omega) \to E_\omega(x)
    \]
    is an isomorphism.
\end{lem}

\begin{proof}
If $f \in H^0(\Sigma;E_\omega)$, then $f \in \mathrm{Ker}(\Delta_\omega)$ which is $1$-dimensional (Lemma \ref{lemm:gap}) and is identified with the span of the constant function $1$, so it does not vanish.
\end{proof}

 We will also need the following lemma.
    \begin{lem}
    \label{lemma:poles}
        Let $x,y\in \Sigma$ be two distinct points. Then there exists a meromorphic $1$-form $\beta_{x,y}$ on $\Sigma$ with exactly 2 simple poles at $x$ and $y$. Moreover, we can choose $\beta_{x,y}$ such that 
        \begin{equation}
            \label{eq:Residue}
    \mathrm{Res}_x(\beta_{x,y})=1,\quad  \mathrm{Res}_y(\beta_{x,y})=-1,
        \end{equation}
        where $\mathrm{Res}_p(\varphi)$ denotes the residue of the form $\varphi$ at $p \in \Sigma$.
    \end{lem}
    \begin{proof}
        Recall that $K$ denotes the canonical line bundle of $\Sigma$; its holomorphic sections are holomorphic $1$-forms. We are interested in sections of $H^0(\Sigma; K(x+y))$, where $x+y$ stands for the effective divisor, and where the notation means that we consider meromorphic $1$-forms with at most simple poles at $x$ and $y$. 
        
        We denote by $\mathcal O$ the trivial holomorphic line bundle, whose holomorphic sections are holomorphic functions. We let $\mathcal O(-x-y)$ denote the space of holomorphic functions vanishing at $x$ and $y$ (this is reduced to $0$). Note that we have $K\otimes \mathcal O(x+y)=K(x+y)$. Applying the Riemann--Roch theorem to the holomorphic line bundle $L=K(x+y)$, we obtain
       \begin{align*}
           \mathrm{dim}(H^0(\Sigma;L))-\mathrm{dim}(H^0(\Sigma;K\otimes L^{-1}))=\mathrm{deg}(L)-\textsl{g}+1,
       \end{align*}
       that is, 
       \[
       \begin{split}
       \mathrm{dim}(H^0(\Sigma;K(x+y))) &=\mathrm{deg}(K)+2+1-\textsl{g}+\mathrm{dim}(H^0(\Sigma;\mathcal O(-x-y))) \\
       & =2\textsl{g}-2+2+1-\textsl{g}=\textsl{g}+1,
       \end{split}
       \]
       where we used that $\mathrm{deg}(K)=2\textsl{g}-2$ and that holomorphic functions of $\Sigma$ are constant to obtain $\mathrm{dim}(H^0(\Sigma;\mathcal O(-x-y)))=0.$ 
       
       Since $\textsl{g}>0$, this shows the existence of a (non-constant) meromorphic $1$-form with at most simple poles a $x$ and $y$. By the residue theorem, the sum of residues has to vanish and hence the $1$-form must have a pole at both $x$ and $y$. We can ensure \eqref{eq:Residue} by rescaling the $1$-form.
    \end{proof}

\subsubsection{Expressing the value at $0$ using the spectral projector on $1$-forms}

We will now use the fact that the dual line bundle $E_\omega^*$ identifies naturally with $E_{-\omega}$; this is still a holomorphic line bundle over $\Sigma$. Let $x,y\in\Sigma$ and consider $\xi \in E_{\omega}(x)$ (resp. $\mu\in E_{-\omega}(y)$). By Lemma \ref{lemma:evaluation}, $\xi$ extends uniquely to a holomorphic section $\xi_\star\in H^0(\Sigma;E_{\omega})$ (resp. $\mu_\star \in H^0(\Sigma;E_{-\omega}) $) by defining for any $z\in \Sigma$,
\begin{equation}
\label{eq:xi*}
\xi_\star(z):=\mathrm{ev}_z^\omega\circ (\mathrm{ev}_x^\omega)^{-1}(\xi) = \sigma_\omega(x,z)\xi.
\end{equation}
Using Proposition \ref{prop:DR}, for any $s\in \CC$ such that $\mathrm{Re}(s)\gg1$, and any $\eps>0$ small enough,
$$ 
\langle \eta_{g,\omega,x,y}(s)\xi,\mu\rangle=e^{-s\varepsilon}\langle R_1^{\omega}(s)\big(\iota_X[S_x\Sigma]\otimes \pi^*\xi_\star\big),e^{\varepsilon \mathcal L_{X^{-\omega}}}([S_y\Sigma]\otimes \pi^*\mu_\star) \rangle+\delta_{x=y} \langle \xi,\mu\rangle.
$$
We emphasize that the previous pairing is $\CC$-linear with respect to both arguments. By Proposition \ref{prop:mero}, the pairing admits a meromorphic extension to $\CC$ with poles included in $\mathrm{Res}^1(\omega)$. Contrary to what we have done in Proposition \ref{prop:primitive}, we cannot deduce the fact that the Poincaré series has a well-defined value at zero from the absence of Jordan blocks at zero. We will use the results of Proposition \ref{prop:gen-res-states} instead.

\begin{prop}
    \label{prop:theo2}
    Let $(\Sigma,g)$ be a hyperbolic surface. There exists a neighborhood $U$ of $0$ in $H^0(\Sigma; K)$ such that for any $x,y\in \Sigma$, $\xi \in E_\omega(x), \mu\in E_{-\omega}(y)$, $\omega \in U\setminus\{0\}$ and $\eps>0$ small enough: 
    \begin{equation}
\label{eq:counterex}
\langle \eta_{g,{\omega},x,y}(0)\xi,\mu\rangle=-\langle \Pi_1^{{\omega}}(0)F_{x,\xi},e^{\varepsilon \mathcal L_{X^{-{\omega}}}}([S_y\Sigma]\otimes \pi^*\mu_\star) \rangle, \end{equation}
where $\Pi_1^{{\omega}}(0)$ denotes the spectral projector defined in \eqref{eq:Pi_k} and $F_{x,\xi}$ is the primitive constructed in Lemma \ref{lem:exactpt}. 
\end{prop}

\begin{rem}
    Note that \eqref{eq:counterex} gives the same result for $x=y$ and $x\neq y$. This fact will be used in the proof of Theorem \ref{theo:2}. 
\end{rem}
\begin{proof}
    The first part of the proof of Proposition \ref{prop:primitive} applies, and in particular, \eqref{eq:1stpart} holds. For $s\in \CC$, we have
    \begin{equation}
    \begin{split}
    \label{eq:name}
        \langle &\eta_{g,\omega,x,y}(s)\xi,\mu\rangle-\delta_{x=y}\langle \xi,\mu\rangle
        \\&=e^{-s\varepsilon}\langle \big([S_x\Sigma]\otimes  \pi^*\xi_\star\big),e^{\varepsilon \mathcal L_{X^{-\omega}}}F_{y,\mu}\rangle-se^{-s\varepsilon}\langle R_2^\omega(s)\big([S_x\Sigma]\otimes  \pi^*\xi_\star\big),e^{\varepsilon \mathcal L_{X^{-\omega}}}F_{y,\mu}\rangle.
        \end{split}
    \end{equation}
    For the first term on the right-hand side of \eqref{eq:name}, we have 
    \begin{equation}
    \label{eq:miracle}
    \begin{split}\lim_{s\to 0}e^{-s\varepsilon}\langle \big([S_x\Sigma]\otimes  \pi^*\xi_\star\big),e^{\varepsilon \mathcal L_{X^{-\omega}}}F_{y,\mu}\rangle= \langle \big([S_x\Sigma]\otimes  \pi^*\xi_\star\big),e^{\varepsilon \mathcal L_{X^{-\omega}}}F_{y,\mu}\rangle=-\delta_{x=y}\langle \xi,\mu \rangle,
        \end{split}
    \end{equation}
    where the last equality follows from (the proofs) of Propositions \ref{prop:xneqy} and \ref{prop:xeqy}. For the second term on the right-hand side of \eqref{eq:name}, we note that, using Lemma \ref{lem:exactpt},
    $$ R_2^\omega(s)\big([S_x\Sigma]\otimes  \pi^*\xi_\star\big)=d^{\nabla_\omega}R_1^\omega(s)F_{x,\xi}.$$
    Proposition \ref{prop:gen-res-states} yields
    \begin{equation}
        \label{eq:structureRes1}
        (\mathcal L_{X^\omega})^2\Pi_1^\omega(0)=0,\quad d^{\nabla_\omega}\mathcal L_{X^\omega}\Pi_1^\omega(0)=0.
    \end{equation}
    As a consequence, we have the Laurent expansion:
    $$R_1^\omega(s)=-\mathcal L_{X^\omega}\Pi_1^\omega(0)/s^2+\Pi_1^\omega(0)/s+R_{\mathrm{hol}}(s), $$
    where $R_{\mathrm{hol}}(s)$ is holomorphic near $s=0.$ It follows from \eqref{eq:structureRes1} that
    $$d^{\nabla_\omega}R_1^\omega(s)F_{x,\xi}=d^{\nabla_\omega}\Pi_1^\omega(0)F_{x,\xi}/s+O(1). $$
    Plugging this into the second term on the right-hand side of \eqref{eq:name}, we find
    \begin{align*}
        se^{-s\varepsilon}\langle R_2^\omega(s)\big([S_x\Sigma]\otimes  \pi^*\xi_\star\big),e^{\varepsilon \mathcal L_{X^{-\omega}}}F_{y,\mu}\rangle&=e^{-s\varepsilon}\langle d^{\nabla_\omega}\Pi_1^\omega(0)F_{x,\xi},e^{\varepsilon \mathcal L_{X^{-\omega}}}F_{y,\mu}\rangle+O(s)
        \\&=e^{-s\varepsilon}\langle \Pi_1^\omega(0)F_{x,\xi},e^{\varepsilon \mathcal L_{X^{-\omega}}}([S_y\Sigma]\otimes \pi^*\mu_\star )\rangle+O(s),
\end{align*}
where we integrated by parts in the second line. Taking the limit when $s\to 0$, inserting the result in \eqref{eq:name}, and using \eqref{eq:miracle}, finishes the proof of the proposition.
\end{proof}

\subsubsection{End of the proof} We now prove Theorem \ref{theo:2}.
\begin{proof}[Proof of Theorem \ref{theo:2}]
   Recall the notations from Lemma \ref{lemma:poles}. Observe that
    \begin{equation}
        \label{eq:ThibaultMagic}
d\beta_{x,z}=\bar{\partial}\beta_{x,z}=2\pi i([x]-[z]),
    \end{equation}
    where $[x] \in \mathcal{D}'(\Sigma,\Lambda^2T^*\Sigma)$ denotes the distributional $2$-form such that $\langle [x],f\rangle = f(x)$ for all functions $f$. In the first equality, we used that $\partial \beta_{x,z}$ has type $(2,0)$ so it vanishes on a surface. Since $M = S\Sigma$ is oriented using $\alpha\wedge d\alpha$, a direct computation shows that $\pi^*([x])=[S_x\Sigma]$, from which we deduce
\begin{equation}
\label{eq:noidea}d(\pi^*\beta_{x,z})=\pi^*d\beta_{x,z}=2\pi i([S_x\Sigma]-[S_z\Sigma]). 
\end{equation}
Recall that $\xi_\star(z)$ is defined in \eqref{eq:xi*}. 
By Leibniz' rule, it holds that:
\begin{align*}d^{\nabla_\omega}(\pi^*\beta_{x,z}\otimes \pi^*\xi_\star)&=d(\pi^*\beta_{x,z})\otimes \pi^*\xi_\star-\pi^*\beta_{x,z}\wedge \nabla_\omega\pi^*\xi_\star
\\&=d(\pi^*\beta_{x,z})\otimes \pi^*\xi_\star-\pi^*(\beta_{x,z}\wedge (\omega\otimes \xi_\star))
\\&=d(\pi^*\beta_{x,z})\otimes \pi^*\xi_\star,
\end{align*}
    where we used in the second equality that $\nabla_\omega\pi^*\xi_\star= \pi^* \nabla_\omega \xi_\star = \pi^*(\omega\otimes \xi_\star)$ and in the last equality that both $\pi^*\beta_{x,z}$ and $\pi^*(\omega\otimes \xi_\star)$ are of type $(1,0)$ so their wedge product vanishes. Together with \eqref{eq:noidea}, we thus obtain
    \begin{equation}
        \label{eq:finalprimitive}
    d^{\nabla_\omega}\pi^*(\beta_{x,z}\otimes \xi_\star)=2\pi i([S_x\Sigma]-[S_z\Sigma])\otimes \pi^*\xi_\star.    
    \end{equation}
    
    \begin{lem}
    \label{lemm:indep}
        Let $x,y, z\in \Sigma$. We have
        \begin{equation}
        \label{eq:magic}\eta_{g,\omega,x,y}(0)=\eta_{g,\omega,z,y}(0)\circ \sigma_\omega(x,z). 
        \end{equation}
    \end{lem}
    
    \begin{proof}
        The statement of the lemma is trivial if $x=z$ so we can suppose that $x\neq z$. Using $\beta_{z,x}=-\beta_{x,z}$, it follows from Proposition \ref{prop:theo2} that
        \[
        \begin{split}\langle \eta_{g,{\omega},x,y}(0)\xi,\mu\rangle-\langle \eta_{g_0,{\omega},z,y}(0)\xi_\star(z),\mu\rangle & =\langle \Pi_1^{{\omega}}(0)(F_{z,\xi}-F_{x,\xi}),e^{\varepsilon \mathcal L_{X^{-{\omega}}}}([S_y\Sigma]\otimes \pi^* \mu_\star) \rangle \\
        & \overset{\eqref{eq:finalprimitive}}{=} \frac{1}{2i\pi}\langle \Pi_1^{{\omega}}(0)\pi^*(\beta_{z,x}\otimes\xi_\star),e^{\varepsilon \mathcal L_{X^{-{\omega}}}}([S_y\Sigma]\otimes \pi^*\mu_\star) \rangle \\
        & = \frac{1}{2i\pi}\langle \beta_{z,x}\otimes \xi_\star,\pi_*\Pi_2^{*,\omega}(0)\big(e^{\varepsilon \mathcal L_{X^{-{\omega}}}}([S_y\Sigma]\otimes \pi^*\mu_\star)\big) \rangle.
        \end{split}
        \]
        In the second equality, the pairing is well-defined because $\mathrm{WF}\big(\pi^*(\beta_{z,x}\otimes \xi_\star)\big)\subset (\mathbb H^* \oplus \mathbb R\alpha)\setminus \{0\}$, and is thus disjoint from $E_u^*$  and $E_s^*$ by \eqref{eq:Vtransverse}. (The fact that the pairing does not depend on the choice of primitive follows from the exactness of  $e^{\varepsilon \mathcal L_{X^{-{\omega}}}}([S_y\Sigma]\otimes \pi^*\mu_\star )$ provided by Lemma \ref{lem:exactpt}, and an argument similar to the one in the proof of Proposition \ref{prop:xneqy}.) The third equality follows from an integration by parts, and $\Pi_2^{*,\omega}(0)$ denotes the spectral projector onto coresonant states, that is, the resonant states on $2$-forms corresponding to $(-X,-\omega)$.
        
        We claim that 
        \begin{equation}
            \label{eq:claim}
    \pi_*\Pi_2^{*,\omega}\big(e^{\varepsilon \mathcal L_{X^{-{\omega}}}}([S_y\Sigma]\otimes \pi^*\mu_\star)\big)\in \mathrm{Span}(\omega\otimes \mu_\star).
        \end{equation}
        Suppose that the claim holds. Then, inserting this identity in the above equality, we find that there is $c\in \CC$ such that 
        $$\langle \eta_{g,{\omega},x,y}(0)\xi,\mu\rangle-\langle \eta_{g,{\omega},z,y}(0)\xi_\star(z),\mu\rangle=c\langle \beta_{z,x}\otimes \xi_\star,\omega \otimes \mu_\star\rangle=0, $$
        where we used that both $\beta_{z,x}$ and $\omega$ are of type $(1,0)$ so their wedge product vanishes. This concludes the proof of the lemma. 
       \end{proof}
        
\begin{proof}[Proof of \eqref{eq:claim}.] We will use the notation  $v:=\Pi_2^{*,\omega}\big(e^{\varepsilon \mathcal L_{X^{-{\omega}}}}([S_y\Sigma]\otimes \pi^*\mu_\star )\big)$. Let us recall that 
        \begin{equation}
    \label{eq:Res2}
    \begin{split}\mathrm{Res}^{*,2,\infty}(\rho_{-\omega},0)&=\mathrm{Res}^{*,2,1}_0(\rho_{-\omega},0)\oplus( \mathrm{Res}^{*,1,2}_0(\rho_{-\omega},0)\wedge \alpha)
    \\&=(\mathrm{Res}^{*,0,1}_0(\rho_{-\omega},0)\wedge d\alpha)\oplus (\mathrm{Res}^{*,1,2}_0(\rho_{-\omega},0)\wedge \alpha),
    \end{split}
        \end{equation}
which follows from \eqref{eq:relalpha} and Proposition \ref{prop:gen-res-states}. Note that because $\alpha$ restricts to $0$ on each circle fiber, one has $\alpha \wedge [S_y\Sigma]=0$, from which we deduce that 
$\alpha \wedge v=0. $
This means that $v\in \mathrm{Res}^{*,1,2}_0(0)\wedge \alpha$ in the decomposition \eqref{eq:Res2}. Next, we recall that by Lemma \ref{lem:exactpt}, $e^{\varepsilon \mathcal L_{X^{-{\omega}}}}([S_y\Sigma]\otimes \pi^*\mu_\star )$ is $d^{\nabla_{-\omega}}$-exact. By Proposition \ref{prop:gen-res-states}, the resonant $1$-forms are $d^{\nabla_{\omega}}$-closed so $v$ pairs trivially with any genuine resonant $1$-state. This means that 
$v\in \mathcal L_{X^{-\omega}}(\mathrm{Res}_0^{*,1,2}(\rho_{-\omega},0))\wedge \alpha.$ Let $h^*$ be such that it spans $\mathrm{Res}_0^{*,0,1}(\rho_{-\omega},0).$ We claim that
$$\mathcal L_{X^{-\omega}}(\mathrm{Res}_0^{*,1,2}(\rho_{-\omega},0))=\mathrm{Span}(d^{\nabla_{-\omega}}h^*).$$ 
 Indeed, we noticed that $d^{\nabla_{-\omega}}h^*\in \mathrm{Res}_0^{*,1,1}(\rho_{-\omega},0)$ pairs trivially with resonant $1$-states since they are all closed by Proposition \ref{prop:gen-res-states}. The orthogonal (with respect to the pairing between resonant and co-resonant states) of $\mathrm{Res}_0^{1,1}(\rho_\omega,0)$ is given by $\mathrm{Ran}(\mathcal{L}_{X^{-\omega}})$ by duality, and this proves that $d^{\nabla_{-\omega}}h^*\in \mathcal L_{X^{-\omega}}(\mathrm{Res}_0^{*,1,2}(\rho_{-\omega},0))$. 

In order to show that it is non-zero, we compute its first Fourier coefficients. Since the connection is pullbacked from $\Sigma,$ we see that $d^{\nabla_{-\omega}}h^*(V)=V(h^*)$ and in particular, from Lemma \ref{lem:renormalization} (applied to $-X$ and $-\omega$), we have
$$(d^{\nabla_{-\omega}}h^*)_{-1}=0,\quad (d^{\nabla_{-\omega}}h^*)_{1}=2i\omega\neq 0. $$
We can now finish the proof of \eqref{eq:claim}. There is $c\in \CC$ such that for any $x\in \Sigma$, 
\begin{align*}
    \pi_*\Pi_2^{*,\omega}\big(e^{\varepsilon \mathcal L_{X^{-{\omega}}}}([S_y\Sigma]\otimes \pi^* \mu_\star)\big)(x)&=-c\int_{S_x\Sigma} V(d^{\nabla_{-\omega}}h^*)(x,\theta)d\theta\otimes \mu_\star
    \\&=-4i\pi c (\omega \otimes \mu_\star),
\end{align*}
where we used the Fourier modes of $d^{\nabla_{-\omega}}h^*$ obtained above. This finishes the proof of \eqref{eq:claim}.
    \end{proof}

    We now conclude the proof of Theorem \ref{theo:2}. Fix $x\in \Sigma$ and $\xi \in E_\omega(x)$. 

    \begin{lem} The section $\Sigma  \ni y\mapsto \eta_{g,\omega,x,y}(0)\xi \in E_\omega(y)$ is a holomorphic section of $E_\omega$.
    \end{lem}
    
    \begin{proof} Let $t\in H^0(\Sigma;E_{-\omega})$ be a nowhere vanishing section; then for all $x,y\in\Sigma$,
    $$\langle \eta_{g,\omega,x,y}(0)\xi,t(y)\rangle=\langle \xi, \eta_{g,-\omega,y,x}(0)t(y)\rangle, $$
    where the pairing is $\CC$-linear with respect to both arguments. Using Lemma \ref{lemm:indep} (applied to $-X$ and $-\omega$), we see that for all $x,y,z\in \Sigma$
    $$ \eta_{g,-\omega,y,x}(0)t(y)=\eta_{g,-\omega,z,x}(0)\circ \mathrm{ev}_z^{-\omega}\circ (\mathrm{ev}_y^{-\omega})^{-1}(t(y))=\eta_{g,-\omega,z,x}(0)t(z), $$
    which shows that the pairing does \emph{not} depend on $y$. This proves that $ y\mapsto \eta_{g_0,\omega,x,y}(0)\xi$ is holomorphic. 
    \end{proof}

    It then follows from Lemma \ref{lemma:evaluation} that $\eta_{g,\omega,x,y} = c_x \sigma_\omega(x,y)$ for some $c_x \in \CC$. Applying Lemma \ref{lemm:indep} one more time with $X$ and $\omega$, we then find that $c_x=c_z$ for $x,z \in \Sigma$, that is,
\begin{equation}
    \label{eq:almsotfinished}
    \exists C(\omega) \in \CC,\quad \forall x,y\in \Sigma,\quad \eta_{g,\omega,x,y}(0)=C(\omega) \sigma_\omega(x,y).
\end{equation}

\begin{lem}
It holds: $C(\omega)=\chi(\Sigma)^{-1}$.
\end{lem}

\begin{proof}
    We adapt the computation of Paternain \cite{Pat00} in the twisted case and compute the integral of the Poincaré series over $\Sigma\times \Sigma.$ We choose two nowhere vanishing holomorphic sections $\xi_\star $ (resp. $\mu_\star$) of $E_\omega$ (resp. $E_{-\omega}$) normalized so that $\langle \xi_\star,\mu_\star\rangle\equiv 1.$ For any $x,y\in \Sigma$ and $\gamma \in \mathcal P_{x,y}$, we have 
    \begin{equation}
    \label{eq:ppp}\langle \tau_\omega(\gamma)\xi_\star,\mu_\star\rangle=\exp\left(-\int_\gamma \omega\right). \end{equation}
    For $s\in \CC$ such that $\mathrm{Re}(s)\gg 1$, let us define
$$I_\omega(s):=\int_{\Sigma}\int_{\Sigma}\langle \eta_{g,\omega,x,y}(s)\xi_\star(x),\mu_\star(y)\rangle d\mathrm{vol}_{g}(x)d\mathrm{vol}_{g}(y). $$
By a direct adaptation of \cite[Proposition 3.2 and Remark 3.4]{Pat00} and \eqref{eq:ppp}, we have
$$I_\omega(s)=\int_{\Sigma}\int_{0}^{+\infty}e^{-st}\sinh(t)\left( \int_{S_x\Sigma}\exp\left(-\int_{(\gamma_{x,v})|_{[0,t]} }\omega\right)dS_x(v)\right)dtd\mathrm{vol}_{g}(x), $$
where $\gamma_{x,v}$ denotes the $g_0$-geodesic defined by $(x,v).$ The inner integral equals
$$\forall x\in \Sigma, \ \forall t\in \mathbb R, \quad \int_{S_x\Sigma}\exp\left(-\int_{(\gamma_{x,v})|_{[0,t]} }\omega\right)dS_x(v)=2\pi. $$
Indeed, if we lift everything to the universal cover, then the integrand is a holomorphic function of the endpoint. By the mean value theorem, the average of the function on the geodesic sphere is equal to $1,$ its value at the origin. 

This means that 
$$I_\omega(s)=2\pi \mathrm{Area}_{g}(\Sigma)\int_0^{+\infty} e^{-st}\sinh(t)dt=\frac{2\pi \mathrm{Area}_{g}(\Sigma)}{s^2-1}. $$
The right-hand side has a meromorphic extension whose value at $s=0$ equals $-2\pi \mathrm{Area}_{g_0}(\Sigma)$. Integrating \eqref{eq:almsotfinished}, we thus obtain, using Gauss--Bonnet, 
\begin{equation*}
    C(\omega) \mathrm{Area}_{g}(\Sigma)^2=- 2\pi \mathrm{Area}_{g}(\Sigma)\ \Longrightarrow \ C(\omega)=-\frac{2\pi}{\mathrm{Area}_{g}(\Sigma)}=\frac{1}{\chi(\Sigma)}.\qedhere
\end{equation*}
\end{proof}
This concludes the proof of the theorem.
\end{proof}

\bibliographystyle{alpha}
\bibliography{biblio}

\end{document}